\documentclass[preprint]{elsarticle}

\usepackage{amsfonts}
\usepackage{amsmath}
\usepackage{graphicx}      
\usepackage{float} 
\usepackage{psfrag} 
\usepackage{subfig}
\usepackage{multicol}
\usepackage{stfloats}
\usepackage{mathtools}
\usepackage{color}

\usepackage{algorithm}
\usepackage{algpseudocode}

\newtheorem{remark}{Remark}[section]

\makeatletter
\newenvironment{method}[1][htb]{%
    \renewcommand{\ALG@name}{Method}
   \begin{algorithm}[#1]%
  }{\end{algorithm}}

\makeatother

\journal{Computer Methods in Applied Mechanics and Engineering}

\begin{document}

\begin{frontmatter}

\title{Sparsity Preserving Optimal Control of Discretized PDE Systems}


\author[mymainaddress]{Aleksandar Haber\corref{mycorrespondingauthor}}
\cortext[mycorrespondingauthor]{Corresponding author}
\ead{aleksandar.haber@gmail.com; aleksandar.haber@csi.cuny.edu}

\author[mysecondaryaddress]{Michel Verhaegen}

\address[mymainaddress]{Department of Engineering Science and Physics, City University of New York, College of Staten Island, New York 10314, USA}
\address[mysecondaryaddress]{Delft Center for Systems and Control, Delft University of Technology, 2628 CD Delft, The Netherlands}

\begin{abstract}
We focus on the problem of optimal control of large-scale systems whose models are obtained by discretization of partial differential equations using the Finite Element (FE) or Finite Difference (FD) methods. The motivation for studying this pressing problem originates from the fact that the classical numerical tools used to solve low-dimensional optimal control problems are computationally infeasible for large-scale systems. Furthermore, although the matrices of large-scale FE or FD models are usually sparse banded or highly structured, the optimal control solution computed using the classical methods is dense and unstructured. Consequently, it is not suitable for efficient centralized and distributed real-time implementations. We show that the \textit{a priori} (sparsity) patterns of the exact solutions of the generalized Lyapunov equations for FE and FD models are banded matrices. The \textit{a priori} pattern predicts the dominant non-zero entries of the exact solution. We furthermore show that for well-conditioned problems, the \textit{a priori} patterns are not only banded but also sparse matrices. On the basis of these results, we develop two computationally efficient methods for computing sparse approximate solutions of generalized Lyapunov equations. Using these two methods and the inexact Newton method, we show that the solution of the generalized Riccati equation can be approximated by a banded matrix. This enables us to  develop a novel computationally efficient optimal control approach that is able to preserve the sparsity of the control law. We perform extensive numerical experiments that demonstrate the effectiveness of our approach.
\end{abstract}


\end{frontmatter}


\section{Introduction}
Large-scale systems with the dynamics described by Partial Differential Equations (PDEs) are ubiquitous in engineering and science. This comes from the fact that the majority of the fundamental physical laws and processes are mathematically described by PDEs. On the other hand, the problem of optimal control of large-scale systems is a long-standing and pressing problem whose solution is important for efficient, safe, and low-cost operation of various mechanical, electrical and physical systems~\cite{zhou1996robust,benner2016low,motter2013,cornelius2013realistic,haber2017state,lee2011,favennec2010use,haber2012}.

The focus of this paper is on the systems whose dynamics can be represented in the \textit{descriptor} state-space form \cite{kunkel1997linear,duan2010analysis,verhaegen1986reduced}: $E\dot{\mathbf{x}}=A\mathbf{x}+B\mathbf{u}$, where $\mathbf{x}\in \mathbb{R}^{n}$ is the state vector, $\mathbf{u}\in \mathbb{R}^{m}$ is the control input vector, $E$, $A\in \mathbb{R}^{n \times n}$, and $B\in \mathbb{R}^{n \times m}$ are the system matrices. The models where $E=I$ are referred to as \textit{non-descriptor} state-space models. In this paper, we mainly focus on the descriptor state-space models resulting from the discretization of linear PDEs using the classical Finite Element (FE) or Finite Difference (FD) methods. Under the term the classical FE method, we understand a Galerkin FE method with locally supported basis functions and with structured or unstructured meshes composed of the triangular or rectangular 2D elements (tetrahedral and hexahedral 3D elements). FE models of various classes of linear PDEs can be written in the descriptor state-space form. For example, the general form of FE models used for the dynamical analysis is~\cite{bathe2006finite}: $M_{1}\ddot{\mathbf{s}}+M_{2}\dot{\mathbf{s}}+M_{3}\mathbf{s}=B_{1}\mathbf{u}$, where $\mathbf{s},\dot{\mathbf{s}}$, and $\ddot{\mathbf{s}}$ are the displacement, velocity, and acceleration vectors, and  $M_{1},M_{2}$, and $M_{3}$ are the mass, damping and stiffness matrices, respectively. The term $B_{1}\mathbf{u}$ is the vector of applied loads. The corresponding descriptor state-space model is 
\begin{align}
\begin{bmatrix}I & 0 \\  0 & M_{1} \end{bmatrix}\dot{\mathbf{x}}=\begin{bmatrix}0 & I  \\ -M_{3} & -M_{2} \end{bmatrix}\mathbf{x}+\begin{bmatrix}0 \\ B_{1} \end{bmatrix}\mathbf{u}, \;\; \mathbf{x}=\begin{bmatrix}\mathbf{s}\\ \dot{\mathbf{s}}  \end{bmatrix}.
\label{descriptorStateSpaceModel}
\end{align}  
Due to the local support of the basis functions, and taking into consideration the fact that the node degrees in the classical FE methods are relatively small,  the FE element matrices $M_{i}$ are usually sparse~\cite{bathe2006finite,rao2017finite}. Furthermore, after suitable permutations (relabeling of the nodes), they usually have a banded structure~\cite{liu1976comparative}. This furthermore implies that the matrices $A$ and $E$ of the descriptor state-space model are usually sparse banded matrices or its partitions are sparse banded matrices. For example, in the case of \eqref{descriptorStateSpaceModel}, the non-zero blocks of the matrices $E$ and $A$ are banded matrices after suitable permutations. For clarity and brevity, the methods proposed in this paper are developed under the assumption that the matrices $A$ are $E$ are sparse banded matrices (as in the case of the discretized heat equation described below), and as we explain later on, the developed methods can be easily generalized to the case where the block partitions of $A$ and $E$ are sparse banded matrices (as in the case of \eqref{descriptorStateSpaceModel}). 

\begin{figure}[H]
\centering
 \includegraphics[scale=0.28,trim=0mm 0mm 0mm 0mm ,clip=true]{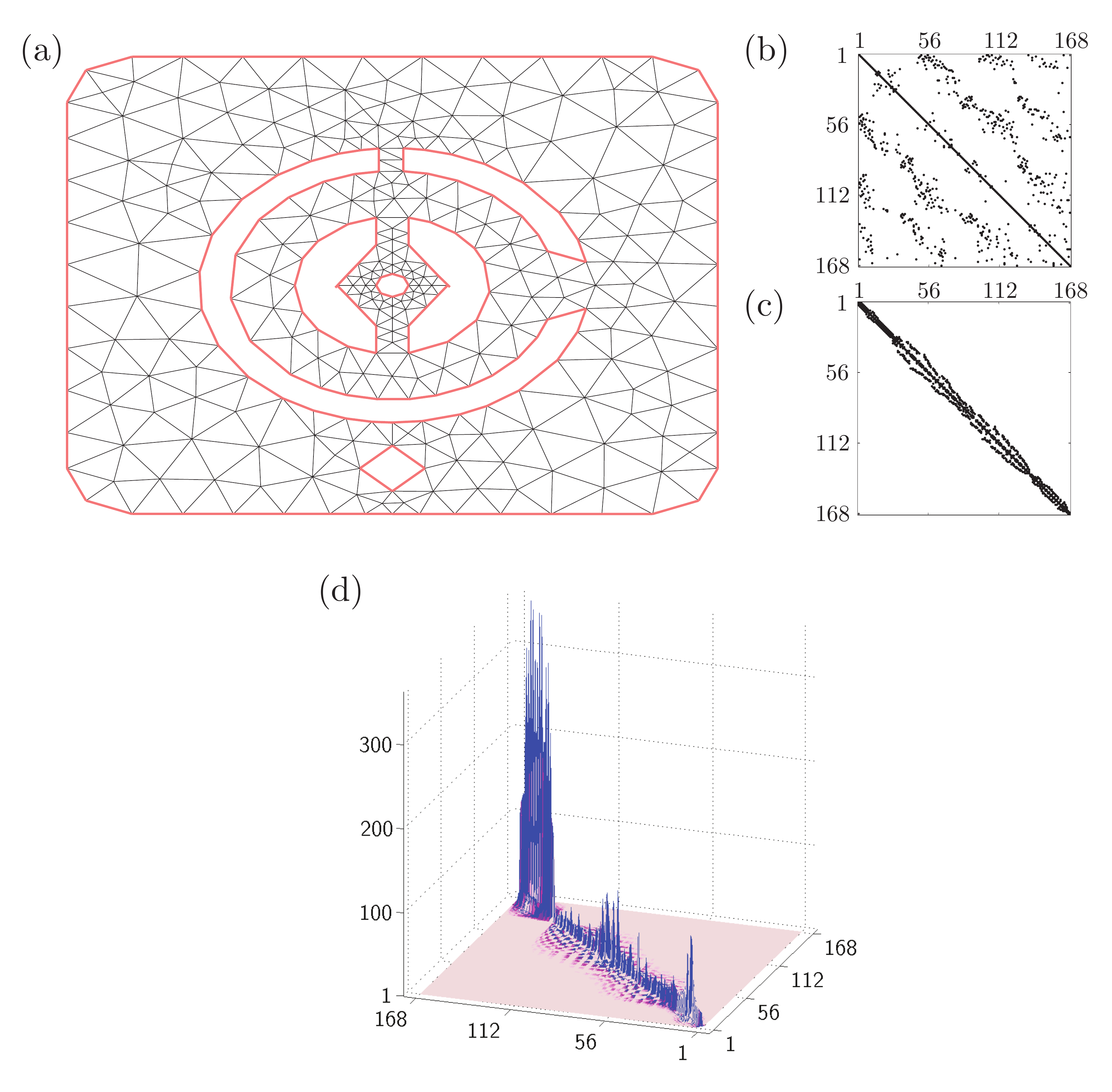}
\caption{(a) A 2D domain and a FE mesh for the discretization of the heat equation. Patterns of (b) the matrix $A$ and (c) the matrix obtained by permuting A using the reverse Cuthill-McKee algorithm. (d) A surface plot of the absolute values of the entries of the solution of the Ricatti equation.}
\label{fig:Gometry}
\end{figure}

To illustrate the structure of the FE matrices, we generate a complex geometry 2D domain shown in Fig.~\ref{fig:Gometry}(a). 
Such a domain can for example represent a cross-section of a building or of an arbitrary mechanical element. We assume that the heat actuators are distributed across the domain, and that we are able to observe the temperature change at predefined points. We assume that the temperature dynamics of such an object is described by the heat equation with Dirichlet boundary conditions and we introduce the FE mesh shown in Fig.~\ref{fig:Gometry}(a). We use tent-shaped linear basis functions. The pattern of the matrix $A$ (the matrix $E$ has a similar pattern) of the descriptor FE state-space model is shown in Fig.~\ref{fig:Gometry}(b). In Fig.~\ref{fig:Gometry}(c) we show the pattern of the matrix $A$, after applying the reverse Cuthill-McKee permutation algorithm~\cite{liu1976comparative}. Similar banded matrix structures appear when the domain is three-dimensional, as well as in the models resulting from the FD methods~\cite{zienkiewicz2005finite}.

The classical numerical algorithms~\cite{anderson2007optimal,kwakernaak1972linear} for solving optimal control problems have $O(n^3)$ computational and $O(n^2)$ memory complexities~\cite{simoncini2016computational,gajic2008lyapunov,kodra2017finding}, and consequently, they are computationally infeasible for large-scale systems. Furthermore, these algorithms produce centralized control solutions which require all sensor data collected from the system to be transferred to a centralized computing unit that calculates and sends control signals to system's actuators. Such an approach usually fails because in the case of a large number of actuators and sensors, an enormous amount of information needs to be transferred and processed in real-time. This implies that in the case of large-scale systems, control has to be performed in a distributed/decentralized manner requiring only local computation and communication between sensors and actuators~\cite{bamieh02distributedcontrol,dandrea2003,jovbamTAC05platoons,linfarjovTAC12platoons,khan2008distributing}. This further implies that the matrices of the optimal controllers need to be sparse and highly structured, mathematically reflecting the need for distributed communication and computation. 

The problem of designing sparse control matrices has been considered by several authors, see for example \cite{linfarjovTAC13admm,linfarjovTAC11al,schuler2011} and follow-up approaches. Due to their $O(n^3)$ computational and $O(n^2)$ memory complexities such methods are computationally infeasible for large-scale systems. Several other methods with lower computational and memory complexities have been developed \cite{justin2009,massioni2009}. However, these methods are applicable only to a narrow class of systems, such as systems composed of strings of subsystems \cite{justin2009}, or decomposable systems \cite{massioni2009}. On the other hand, it has been observed that in the case where the number of system's inputs and outputs is relatively small compared to $n$, the solution of the algebraic Riccati equation necessary to compute the Linear-Quadratic (LQ) control law, can be factored into a product of low-rank matrices \cite{benner2008numerical,bini2011numerical,benner2004solving}. Such an LQ control law can be computed and implemented relatively efficiently. The applicability of these methods is restricted mainly because large-scale systems usually have a large number of inputs and outputs. Moreover,  the computed LQ feedback matrix is dense and thus it is difficult to implement the resulting controller in distributed or decentralized manners.

The main computational bottleneck in solving the optimal control problems originates from the high complexity of the numerical methods for solving the Lyapunov and Riccati matrix equations~\cite{simoncini2016computational,mehrmann1991autonomous}. In this paper, the Lyapunov equation defined for descriptor state-space models is referred to as \textit{the Generalized Lyapunov (GL) equation}, whereas the one defined for the non-descriptor state-space models (where $E=I$), is referred to as \textit{the Non-Generalized Lyapunov (NGL) equation}. In the  same spirit we refer to Riccati equations as \textit{Generalized Riccati (GR)} and \textit{Non-Generalized  Riccati (NGR) equations}. Recently, it has been shown that the solutions of the NGL equations with symmetric banded matrices belong to the class of localized matrices \cite{Haber2016sparse,simoncini2015lyapunov}. Loosely speaking, localized matrices are characterized by the fact that a small number of matrix entries is significantly larger than other entries, and the dominant entries are usually grouped in localized matrix regions \cite{benzi2015decay,benzi2007,haber2016framework,haber2013moving,simoncini2015lyapunov,pozza2016decay}. Furthermore, the entries are decaying away from these regions. Off-diagonally decaying matrices are typical examples of localized matrices \cite{benzi1999bounds,demko1984,benzi2016survey}.
 An example of an off-diagonally decaying matrix is shown in Fig.~\ref{fig:Gometry}(d). 
 
 The localization and spatially decaying phenomena of the solutions of distributed control problems  for infinite-dimensional systems have been analyzed in \cite{motee2013measuring,motee2009distributed,motee2017sparsity}. It has been shown that the solutions of the operator NGL and NGR equations are spatially decaying. The operator NGL and NGR equations are the infinite-dimensional counterparts of finite-dimensional NGL and NGR equations~\cite{opmeer2013finite,ito1998approximation,grubivsic2014eigenvalue}. The extension of these approaches to finite-dimensional systems and to GL and GR equations, is not straightforward. Furthermore, these approaches do not take into consideration the influence of the condition numbers of the system matrices on the spatially decaying phenomenon. Finally, they do not deal with the pressing problem of efficiently computing the distributed controllers.  An approach for computing low-rank solutions to the operator NGL equations has been developed in~\cite{opmeer2013finite}. However, this approach is not suitable for distributed implementation.
 
 In \cite{Haber2016sparse}, we have shown that the decay rate of the entries of the solution of the NGL equation is faster for smaller condition numbers of $A$. By exploiting this localization phenomena, we have developed two efficient methods for computing approximate sparse banded solutions of the NGL equations. The importance of the results of \cite{Haber2016sparse} lies in the fact that the Riccati equation can be solved by solving a series of Lyapunov equations. It should be expected that the solution of the Riccati equation is localized and that it can be approximated by a banded matrix in a computationally efficient manner. Consequently, it might be possible to efficiently approximate the solution of the LQ optimal control problem by a banded matrix. However, the results and the methods presented in \cite{Haber2016sparse} are not applicable either to state-space models with non-symmetric matrices, or to GL and GR equations. On the other hand, numerical evidence indicates that the solutions of GL and GR equations are localized. As an illustration, Fig.~\ref{fig:Gometry}(d) shows the surface plot of the solution of the GR equation for the discretized heat equation and for the geometry defined in panel (a). This important observation motivates us to develop the methods presented in this paper. 

In this paper, we show that the \textit{a priori} (sparsity) patterns of the exact solutions of the GL equations for FE or FD descriptor state-space models are banded matrices. The \textit{a priori} pattern predicts the dominant non-zero entries of the exact solution~\cite{huckle1999apriori}. We furthermore show that for well-conditioned problems (the notion of a well-conditioned problem is introduced in Section~\ref{sectionAprioriSparsityPattern}), the \textit{a priori} patterns are not only banded but also sparse matrices. On the basis of these results, we develop two computationally efficient methods for computing sparse approximate solutions of GL equations. In contrast to \cite{Haber2016sparse}, the methods presented  in this paper are applicable to descriptor state-space models with non-symmetric matrices. Using the developed methods and the inexact Newton method~\cite{benner2016inexact,feitzinger2009inexact}, we show that the solution of the generalized Riccati equation can be approximated by a banded matrix. This enables us to approximate the feedback matrices of the LQ controllers by (sparse) banded matrices and to develop a novel computationally efficient approach for optimal control of large-scale systems.  We perform extensive numerical experiments to demonstrate the effectiveness of the proposed approach. The MATLAB codes and the models used in our simulations can be found in~\cite{haber2018codes}.

The main focus of this paper is on the development of the algorithms and on the investigation of their numerical performance. In problem formulations that differ from ours, the influence of the approximation errors (originating from approximately solving NGL equations) on the stability and convergence properties of the inexact Newton method, has been studied in~\cite{benner2016inexact,feitzinger2009inexact}. In our future work, we will use these results to analyze the stability and convergence properties of the inexact Newton  method used in this paper.

The paper is organized as follows. In Section \ref{sectionProblemFormulation} we formulate the control problem. In Section \ref{sectionAprioriSparsityPattern} we present a method for determining the \textit{a priori} pattern of the solution of the  Lyapunov equation, which is used in Section \ref{sectionLyapunovSolution} to develop the computationally efficient approximation methods. In Section \ref{sectionNumericalExperiments} we present numerical results, and finally, in Section \ref{sectionConclusions} we present the conclusion and discuss the future work.

\section{Problem Formulation}
\label{sectionProblemFormulation}
\subsection{Notation and Preliminaries}
We define the notation used in this paper. The notation $X=[x_{i,j}]$ denotes a matrix whose $(i,j)$ entry is $x_{i,j}$. The exponential of the matrix $X$ is denoted by $\exp\left(X\right)$. The Frobenius and 2-norms are denoted by $\left\| X\right\|_{F}$ and $\left\|X \right\|_{2}$, respectively. The maximal and minimal eigenvalues are denoted by $\lambda_{\max}\left(X\right)$ and $\lambda_{\min}\left(X\right)$, respectively. Similarly, the maximal and minimal singular values are denoted by $\sigma_{\max}\left(X\right)$ and $\sigma_{\min}\left(X\right)$, respectively. The condition number is denoted by $\kappa\left(X\right)$. Let $Q=[q_{i,j}]\in \mathbb{R}^{n\times n}$ be an arbitrary matrix. The projection operator, projecting the matrix $Q$ onto the pattern of an arbitrary matrix $X=[x_{i,j}]\in \mathbb{R}^{n\times n}$, is defined by
\begin{align}
\mathcal{V}_{X}\big[ Q \big]=\left\{ \begin{array} {rl}
q_{i,j}, &\;\; x_{i,j}\ne 0 \;\;   \\
0, &  \;\;  x_{i,j}=0  \end{array} \right. 
\label{sparsificationOperatordefinition}
\end{align}
where $\mathcal{V}_{X}\big[Q\big]\in \mathbb{R}^{n\times n}$.
\subsection{Control Problem Formulation}
For clarity, the methods of this paper are developed assuming a concrete large-scale optimal control problem in which the system matrices $E$ and $A$ have sparse banded structures (after suitable permutations). The generalization of the developed methods to the descriptor state-space models \eqref{descriptorStateSpaceModel} in which the blocks of the system matrices $E$ and $A$ are sparse banded matrices is explained in Remark~\ref{remarkGenerati} in Section~\ref{sectionLyapunovSolution}.

Despite the fact that the matrices of FE or FD models obtained for example using COMSOL or Ansys modeling software are sparse, they are unstructured, and as such, they might not be suitable for the design of distributed control algorithms. This motivates us to first explain a procedure for transforming unstructured FE models into banded state-space models, more suitable for the efficient control design. Consider a 2D domain with the FE mesh shown in Fig.~\ref{fig:Gometry}(a). We assume that heat actuators and temperature sensors are randomly distributed over the domain. In practice, localized heat actuators based on the resistive heaters can be placed
at predefined positions, and the thermocouples can be used as temperature sensors. The control objective is to keep the temperature distribution constant using the feedback information provided by the sensors. Similar control problems appear in various real-life applications, such as for example in control of thermally actuated deformable mirrors used in adaptive optics applications \cite{haber2014estimation,haber2013identification,saathof2015actuation,haber2013predictive,polo2013linear}.
The temperature dynamics is governed by the 2D heat equation with constant parameters. For simplicity, we assume that the temperature is constant and equal to zero on the boundaries (Dirichlet boundary conditions), for more details see Section~\ref{sectionNumericalExperiments}. The discretized model has the following form: $E_{1}\dot{\mathbf{z}}(t)=A_{1}\mathbf{z}(t)$, where $\mathbf{z}(t)\in \mathbb{R}^{n}$ is the vector of temperatures at the discretization nodes, and $E_{1}, A_{1} \in\mathbb{R}^{n\times n}$ are the system matrices. The pattern of the matrix $A_{1}$ is shown in Fig.~\ref{fig:Gometry}(b) (the matrix $E_{1}$ has a similar pattern). We permute the matrices of the state equation into a banded form. The permutation of $A_{1}$ is defined by $A=PA_{1}P^{T}$, $A\in\mathbb{R}^{n\times n}$, where $P\in \mathbb{R}^{n\times n}$ is a permutation matrix of the reverse Cuthill-McKee algorithm \cite{liu1976comparative}. The banded pattern of $A$ is shown in Fig.~\ref{fig:Gometry}(c). The matrix $P$ defines a coordinate transformation $\mathbf{x}(t)=P\mathbf{z}(t)$, where $\mathbf{x}(t)\in \mathbb{R}^{n\times n}$ is the reordered state vector. Taking into account that permutation matrices are orthogonal ($P^{T}P=PP^{T}=I$), we have $\mathbf{z}(t)=P^{T}\mathbf{x}(t)$, and consequently, $PE_{1}P^{T}\dot{\mathbf{x}}(t)=PA_{1}P^{T}\mathbf{x}(t)$. The last equation can be written compactly $E\dot{\mathbf{x}}(t)=A\mathbf{x}(t)$, where $E=PE_{1}P^{T}$, and $E\in \mathbb{R}^{n\times n}$. The pattern of $E$ is similar to the pattern of $A$. We modify the state equation by introducing the control input and add an output equation. The resulting state-space model is 
\begin{align}
E\dot{\mathbf{x}}(t)=A\mathbf{x}(t)+B\mathbf{u}(t), \;\;\; \mathbf{y}(t)&=C\mathbf{x}(t), \notag
\end{align}
where $\mathbf{u}\in \mathbb{R}^{m}$ is the control input vector, $\mathbf{y}\in \mathbb{R}^{r}$ is the observed output vector, $B\in\mathbb{R}^{n\times m}$ and $C\in\mathbb{R}^{r\times n}$ are the system matrices. Assuming that approximately $50\%$ of the mesh points are directly observed and controlled, the patterns of the matrices $B$ and $C$ are shown in Fig.~\ref{fig:SparsityBC} in Section~\ref{sectionNumericalExperiments}. The methods developed in this paper are applicable to systems with the matrices $B$ and $C$ having similar diagonal-like patterns.

 We are interested in finding the control input vector producing a desired temperature profile. That is, we are interested in the \textit{set point tracking problem} \cite{anderson2007optimal,kwakernaak1972linear}. For presentation clarity, we will use a simple set point tracking problem formulation~\cite{kwakernaak1972linear}. Let $\mathbf{u}_{d}$ be the control input that brings the system to the desired state $\mathbf{x}_{d}$ and to the desired output $\mathbf{y}_{d}$ in the steady state. That is, the input $\mathbf{u}_{d}$ should satisfy $ A\mathbf{x}_{d}+B\mathbf{u}_{d}=0$ and $\mathbf{y}_{d}=C\mathbf{x}_{d}$, either exactly or in the least-squares sense. Next, we introduce the error variables $\tilde{\mathbf{u}}(t)=\mathbf{u}(t)-\mathbf{u}_{d}$,  $\tilde{\mathbf{x}}(t)=\mathbf{x}(t)-\mathbf{x}_{d}$, and  $\tilde{\mathbf{y}}(t)=\mathbf{y}(t)-\mathbf{y}_{d}$. It is easy to verify that the error variables satisfy:
\begin{align}
 E\dot{\tilde{\mathbf{x}}}(t) &=A\tilde{\mathbf{x}}(t)+B\tilde{\mathbf{u}}(t), \label{errorVariablesState} \\
 \tilde{\mathbf{y}}(t)&=C\tilde{\mathbf{x}}(t). \label{errorVariablesOutput}
\end{align}
The LQ control problem is defined by  \cite{mehrmann1991autonomous}:
\begin{align}
&\min_{\tilde{\mathbf{u}}(t)}\int_{0}^{\infty}\tilde{\mathbf{y}}(t)^{T}Q \tilde{\mathbf{y}}(t)+\tilde{\mathbf{u}}(t)^{T}R\tilde{\mathbf{u}}(t)\mathrm{d}t,
\label{theCostFunction} \\
& \text{subject to \eqref{errorVariablesState} and \eqref{errorVariablesOutput}},  \notag
\end{align}
where $Q\in \mathbb{R}^{r\times r}$ and $R\in \mathbb{R}^{m\times m}$ are the weighting matrices. We assume that $E$ is nonsingular, $R$ is positive-definite, $(E,A,B)$ strongly stabilizable, and $(E,A,C^{T}QC)$ is strongly detectable. Furthermore, for simplicity and without the loss of generality, we assume that the matrices $Q$ and $R$ are diagonal. Under these assumptions the solution of the LQ control problem exists, it is unique, and it is given by the feedback law~\cite{mehrmann1991autonomous}:
\begin{align}
\tilde{\mathbf{u}}(t)=-F\tilde{\mathbf{x}}(t),\;\; F=R^{-1}B^{T}ZE,
\label{solutionLQproblemControlLaw}
\end{align}
where $Z$ is the unique positive-definite solution of the GR equation $\mathcal{D}\left[Z\right]=0$, where 
\begin{align}
\mathcal{D}\left[Z\right]=
C^{T}QC + E^{T}ZA + A^{T}ZE - 
E^{T}ZBR^{-1}B^{T}ZE.
\label{riccatiOperator}
\end{align}
Furthermore, the (closed loop) dynamics of the system obtained with this control
\begin{align}
E\dot{\tilde{\mathbf{x}}}(t)=\bar{A}\tilde{\mathbf{x}}(t),\;\; \bar{A}=A-BF,
\label{closedLoopStability}
\end{align}
is asymptotically stable. We solve the Riccati equation $\mathcal{D}\left[Z\right]=0$ using the \textit{inexact Newton method} \cite{benner2016inexact,feitzinger2009inexact}. It can be easily shown that the \textit{Fr\'echet derivative} of $\mathcal{D}\left[Z\right]$ is 
\begin{align}
&\mathcal{D}^{'}_{Z}\left[Y\right]=E^{T}YA+A^{T}YE-E^{T}ZBR^{-1}B^{T}YE-E^{T}YBR^{-1}B^{T}ZE. 
\label{FrechetDerivative}
\end{align}
In the $k$-th iteration, the inexact Newton method solves
\begin{align}
\mathcal{D}^{'}_{\hat{Z}_{k-1}}\left[\hat{Z}_{k}-\hat{Z}_{k-1}\right]=-\mathcal{D}\left[\hat{Z}_{k-1}\right]+S_{k},
\label{newtonStep}
\end{align}
where $\hat{Z}_{k}$ is an unknown matrix, and the matrix $\hat{Z}_{k-1}$ is determined in the previous iteration $k-1$. In \eqref{newtonStep}, $S_{k}$ is a residual matrix, quantifying the "inexactness" of the Newton step. The convergence properties of the inexact Newton methods have been studied in~\cite{feitzinger2009inexact,benner2016inexact,chehab2017inexact}. When $S_{k}=0$, the equation \eqref{newtonStep} becomes the \textit{classical Newton iteration}. After some transformations, the equation \eqref{newtonStep} becomes the GL equation
\begin{align}
E^{T}\hat{Z}_{k}\bar{A}_{k-1}+\bar{A}_{k-1}^{T}\hat{Z}_{k}E=P_{k}+S_{k},
\label{LyapunovFinal}
\end{align}
where 
\begin{align}
& P_{k}=-C^{T}QC-\hat{F}_{k-1}^{T}R\hat{F}_{k-1}, \; \bar{A}_{k-1}=A-B\hat{F}_{k-1},\;\; \hat{F}_{k-1}=R^{-1}B^{T}\hat{Z}_{k-1}E.
\label{LyapunovFinalExplanation}
\end{align}
 If the residual term is zero, then $\hat{Z}_{k}$ becomes $Z_{k}$, that exactly solves
\begin{align}
E^{T}Z_{k}\bar{A}_{k-1}+\bar{A}_{k-1}^{T}Z_{k}E=P_{k}.
\label{lyapunovEquationReWritten}
\end{align}
That is, the GR equation is actually solved by solving the GL equations for every iteration $k$ of the Newton method. The inexact Newton method is summarized in Algorithm~\ref{algorithmNewton}. 

The initial guess in Algorithm~\ref{algorithmNewton} should be chosen as a scaled identity matrix or as a sparse banded matrix. This algorithm ensures that the matrices $\hat{Z}_{k}$ are always banded (a method for selecting the structure and the bandwidth of the approximate solution in the second step of Algorithm~\ref{algorithmNewton} is presented in Section~\ref{sectionAprioriSparsityPattern}). Due to this fact, and taking into consideration that the matrices $A$ and $E$ are banded, and that the matrices $Q$, $R$, and $R^{-1}$ are diagonal, the approximation of the feedback control matrix, denoted by $\hat{F}_{k}$, as well as the matrix $P_{k}$, are also banded matrices. The matrices $C$ and $B$, with the patterns shown in Fig.~\ref{fig:SparsityBC}, as well as the matrices $C$ and $B$ with similar diagonal-like patterns, do not fundamentally change the banded structure of the matrices $\hat{Z}_{k}$, $\hat{F}_{k}$, and $P_{k}$.

 \begin{algorithm}[H]
\caption{Inexact Newton Algorithm}
\textbf{Input:} The matrices $E$, $A$, $B$, $C$, $Q$, and $R$.  \\
\textbf{Output:} Banded approximate solution $\hat{Z}$ of the GR equation.
\begin{enumerate}
\item  Choose an initial guess $\hat{Z}_{0}$.
\item For $k=1,2,\ldots, N_{max}$, compute a banded matrix $\hat{Z}_{k}$ that approximately solves the GL equation (up to a residual $S_{k}$ in \eqref{LyapunovFinal}):  
\begin{align}
E^{T}Z_{k}\bar{A}_{k-1}+\bar{A}_{k-1}^{T}Z_{k}E=P_{k}, \notag
\end{align}
 where $N_{max}$ is the maximal number of iterations. Terminate the loop if the convergence tolerance is satisfied.
 \end{enumerate}
 \label{algorithmNewton}
\end{algorithm}
 
 \textit{The main research question is: Under which conditions does the GL equation admit banded approximate solutions? In addition, can the bandwidth and the structure that capture the most dominant entries of the exact solution to the GL equation be predicted a priori? The main challenge is to develop efficient methods for computing banded approximate solutions to the GL equations.} For example, the numerical evidence shows that the exact solutions of the GL and GR equations for the FE heat equation model belong to the class of localized matrices (Fig.~\ref{fig:Gometry}(d) shows the solution of the GR equation), implying that the GL and GR equations admit banded approximate solutions. However, despite the fact that similar numerical observations can be made for other types of systems~\cite{haber2014estimation}, in order to develop efficient numerical algorithms, we need analysis methods that can provide us with deeper understanding of these phenomena.
 
  In the next section we will present a method for computing the \textit{a priori} pattern of the exact solution $Z_{k}$. The \textit{a priori} pattern predicts the dominant non-zero entries of the exact solution~\cite{huckle1999apriori}. This method enables us to show that the banded approximate solution exists due to the fact that the system matrices are banded. For well-conditioned problems (the notion of a well-conditioned problem is introduced in Section~\ref{sectionAprioriSparsityPattern}), the \textit{a priori} pattern is not only banded but also a sparse matrix, implying that the most dominant entries of the exact solution belong to a sparse banded matrix pattern. These results are used in Section~\ref{sectionLyapunovSolution} to develop computationally efficient methods for approximating the solutions of the GL equations. The developed methods are used in the second step of Algorithm~\ref{algorithmNewton}.

\section{A Priori Patterns}
\label{sectionAprioriSparsityPattern}
In this section we  present a method for computing the \textit{a priori} patterns of the exact solutions to the GL equations. The \textit{a priori} pattern will be used to significantly reduce the dimensionality and the complexity of the approximation problem and to develop the efficient approximation methods in Section~\ref{sectionLyapunovSolution}. These reductions will be accomplished by eliminating or neglecting the entries of the exact solution that do not belong to the \textit{a priori} pattern.

In the analysis presented here, we are mainly interested in the non-zero structure of the \textit{a priori} pattern, but not in the exact numerical values of its entries. Furthermore,  we ignore possible zero entries created by coincidence in the numerical values of the matrix entries. This is a standard assumption that is equivalent to the assumption that all the entries of the matrices used for the \textit{a priori} pattern analysis are positive, see for example~\cite{huckle1999apriori} and references therein (see also Remark~\ref{explodingRemark}).

The fact that the numerical evidence shows that the matrix $Z_{k}$ belongs to the class of localized matrices, might lead to the conclusion that we can simply assume a banded \textit{a priori} pattern of $Z_{k}$, with a completely dense bandwidth region, and proceed with the development of the approximation algorithms.  \textit{However, the pattern bandwidth is difficult to guess. Since we are dealing with large-scale problems, \textit{a priori} patterns with completely dense bandwidth regions and large bandwidths will not necessarily produce significant computational savings. In many cases, a sparser banded matrix would be a much better choice that could lead to significant computational savings.} In \cite{Haber2016sparse} we have developed a method for predicting the \textit{a priori} patterns of approximate solutions of NGL equations. Here, we modify and generalize such a method to GL equations. This method is summarized in Algorithm~\ref{algorithmApriori}. 

 Applying the $\text{vec}\left(\cdot\right)$ operator to \eqref{lyapunovEquationReWritten}, we obtain
\begin{align}
M\mathbf{z}=\mathbf{p},
\label{vectorizedLyapunovEquation}
\end{align}
where $M\in \mathbb{R}^{n^2 \times n^2}$, $\mathbf{z}\in \mathbb{R}^{n^2}$, and $\mathbf{p}\in \mathbb{R}^{n^2}$ are defined as follows
\begin{align}
&M=\bar{A}^{T}_{k-1}\otimes E^{T}+E^{T}\otimes \bar{A}_{k-1}^{T},\; \label{vectorizedLyapunovEquationExplanation} \\ 
&\mathbf{z}=\text{vec}\left(Z_{k}\right),\; \mathbf{p}=\text{vec}\left(P_{k} \right).
\label{vectorizedLyapunovEquationExplanation2}
\end{align}
Due to the fact that the dimension of the system \eqref{vectorizedLyapunovEquation} is extremely large, \textit{the matrix $M$ cannot be formed explicitly. Moreover, it might not be possible to compute the solution $\mathbf{z}=M^{-1}\mathbf{p}$ by using the direct solution methods~\cite{davis2016survey} or by using preconditioned iterative techniques~\cite{Benzi2002survey}}. Despite this obstacle, the fact that $M^{-1}$ can be formally expressed as the sum of powers of $M$ can help us gain more insight into the structure of the \textit{a priori} pattern. Similarly to \cite{Haber2016sparse}, we use the Newton-Schulz iteration to analyze the \textit{a priori} pattern. The Newton-Schulz method approximates $M^{-1}$ using the following iteration $X_{l+1}=X_{l}\left(2I-MX_{l}\right)$, where $l=0,1,2,\ldots$\;, is the iteration index, and $X_{l}$ is an approximation of $M^{-1}$. Let the approximation error of the Newton-Schulz iteration be quantified by $E_{l}=I-MX_{l}$. It can be easily shown that if $X_{0}=\Big( 2/\big(\sigma_{max}^{2}(M)+\sigma_{\min}^{2}(M) \big)\Big)M^{T}$, then 
\begin{align}
\left\|E_{l}\right\|_{2} \le \left(\frac{\kappa\left(M\right)^2-1 }{\kappa\left(M\right)^2+1 } \right)^{2^{l}}.
\label{newtonSchulzIterationError}
\end{align}
 Using the Netwton-Schulz iteration and the aforementioned initial guess, by back substitution it can be shown that the pattern of $X_{l}$ is determined by the pattern of the following matrix
\begin{align}
\underbar{X}_{l}=\left(I+M^{T}M+\ldots+\left(M^{T}M\right)^{w} \right)M^{T},
\label{sparsityPatternX}
\end{align}
where $w=2^{l}-1$. From \eqref{newtonSchulzIterationError} we see that if the matrix $M$ is well-conditioned, then for a relatively small $l$ we obtain a good approximation to $M^{-1}$. This is a direct consequence of the quadratic convergence rate of the Newton-Schulz iteration described by \eqref{newtonSchulzIterationError}. This, together with \eqref{sparsityPatternX}, implies that for a well-conditioned $M$, we can get a relatively good guess of the \textit{a priori} pattern of $M^{-1}$ by summing only a few powers of $M^{T}M$ and by multiplying the resulting sum with $M^{T}$. Since the matrices $A$ and $E$ are sparse banded, the matrix $M$ is sparse (multi) banded. \textit{All this implies that for a well-conditioned $M$, the \textit{a priori} pattern that captures the most dominant entries of $M^{-1}$ is a sparse (multi) banded matrix.} In the general case, the fill-in of the matrix $\underbar{X}_{l}$ will increase with $l$ or consequently, with $w$. However, as long as $w\ll n$, the matrix $\underbar{X}_{l}$ will be a sparse (multi) banded matrix.

In practice, the matrix $M$ cannot be formed explicitly due to its large-scale nature. Consequently, we cannot compute the \textit{a priori} patterns of the vector $\mathbf{z}$ and the matrix $Z_{k}$. This problem can be resolved by reversing the $\text{vec}$ operation used to form \eqref{vectorizedLyapunovEquation}. Since the approximate pattern of $M^{-1}$ is determined by  $\underbar{X}_{l}$, from \eqref{vectorizedLyapunovEquation} and \eqref{sparsityPatternX}, we have
\begin{align}
&\underline{\mathbf{z}}= \mathbf{g}+M^{T}M\mathbf{g}+\ldots+\left(M^{T}M\right)^{w}\mathbf{g},
\label{sparsityPatternAnalysis}
\end{align}
where $\underline{\mathbf{z}}$ is a vector whose non-zero entries determine an \textit{a priori} pattern of $\mathbf{z}$, and $\mathbf{g}=M^{T}\mathbf{p}$. From \eqref{sparsityPatternAnalysis}, we can reconstruct the \textit{a priori} pattern of $Z_{k}$. Namely, by reversing the action of the $\text{vec}\left( \cdot \right)$ operator on \eqref{sparsityPatternAnalysis}, we obtain 
\begin{align}
\underline{Z}=G_{1}+G_{2}+\ldots+G_{w+1},
\label{inverseVecOperator}
\end{align}
where $G_{i}$, $i=1,2,\ldots,w+1$, is the matrix corresponding to $(M^{T}M)^{i-1}\mathbf{g}$ term in \eqref{sparsityPatternAnalysis}, and the matrix $\underline{Z}$ is the \textit{a priori} pattern of $Z_{k}$. The matrix terms $G_{i}$ can be computed recursively as follows. Let us start with the first entry in the sum $\mathbf{g}=M^{T}\mathbf{p}$. From the Kronecker sum structure of $M$ given in \eqref{vectorizedLyapunovEquationExplanation}, we can conclude that the matrix form of this vector term, denoted by $G_{1}$, becomes 
\begin{align}
G_{1}=EP_{k}\bar{A}_{k-1}^{T}+\bar{A}_{k-1}P_{k}E^{T}.
\label{G1}
\end{align}
By using the basic properties of the Kronecker product, it can be shown that the vector $M^{T}M\mathbf{g}$ has the following matrix form:
\begin{align}
&G_{2}=E\left( E^{T}G_{1}\bar{A}_{k-1}+\bar{A}_{k-1}^{T}G_{1}E\right)\bar{A}_{k-1}^{T}+\bar{A}_{k-1}\left( E^{T}G_{1}\bar{A}_{k-1}+\bar{A}_{k-1}^{T}G_{1}E\right)E^{T}.
\label{G2}
\end{align}
Using a similar reasoning, it can be shown that the matrix $G_{i+1}$ corresponding to the vector $\left(M^{T}M\right)^{i}\mathbf{g}$ can be computed recursively (see Remark~\ref{explodingRemark})
\begin{align}
&G_{i+1}=E\left( E^{T}G_{i}\bar{A}_{k-1}+\bar{A}_{k-1}^{T}G_{i}E\right)\bar{A}_{k-1}^{T}+\bar{A}_{k-1}\left( E^{T}G_{i}\bar{A}_{k-1}+\bar{A}_{k-1}^{T}G_{i}E\right)E^{T}.
\label{G23}
\end{align}
Our numerical experience shows that in order to get a good estimate of the \textit{a priori} pattern, an identity matrix has to be added to \eqref{inverseVecOperator}, resulting in
\begin{align}
\underline{Z}_{w}=I+G_{1}+G_{2}+\ldots+G_{w+1}.
\label{finalAprioriSparsityPattern}
\end{align}
Using \eqref{finalAprioriSparsityPattern} we can predict the \textit{a priori} pattern of the matrix $Z_{k}$. The parameter $w$ is a user choice. The procedure for computing the \textit{a priori} pattern in summarized in Algorithm~\ref{algorithmApriori}.
\begin{algorithm}[H]
\caption{Computation of the \textit{a priori} pattern $\underline{Z}_{w}$}
\textbf{Input:} the matrices $\bar{A}_{k-1}$, $E$ and $P_{k}$ and a relatively small positive integer $w$ (see Remark~\ref{remarkFillIn} and Remark~\ref{explodingRemark}). \\
\textbf{Output:} \textit{a priori} pattern $\underline{Z}_{w}$.
\begin{enumerate}
\item  Compute $G_{1}=EP_{k}\bar{A}_{k-1}^{T}+\bar{A}_{k-1}P_{k}E^{T}$.
\item For every $i=1,2,\ldots, w$, compute recursively
\begin{align}
&G_{i+1}=E\left( E^{T}G_{i}\bar{A}_{k-1}+\bar{A}_{k-1}^{T}G_{i}E\right)\bar{A}_{k-1}^{T}\notag \\
&+\bar{A}_{k-1}\left( E^{T}G_{i}\bar{A}_{k-1}+\bar{A}_{k-1}^{T}G_{i}E\right)E^{T}. \notag 
\end{align}
\item Compute $\underline{Z}_{w}=I+G_{1}+G_{2}+\ldots+G_{w+1}.$
\end{enumerate}
\label{algorithmApriori}
\end{algorithm}
For a relatively small value of $w$ ($w\ll n$), the \textit{a priori} pattern can be computed very quickly. Namely, since the matrices $A$ and $E$ are sparse banded, the \textit{a priori} pattern can be computed with the $O(n)$ computational and memory complexities. Consequently, our method can be effectively used for large-scale problems. Due to the fact that the matrices $A$, $E$, and $P_{k}$ are banded, the matrix $\underline{Z}_{w}$ is also banded (for $w\ll n$), however, its fill-in increases by increasing $w$, see also Remark~\ref{remarkFillIn}. On the other hand, we have shown that for a well-conditioned $M$, $w$ is generally small, implying that the \textit{a priori} pattern of $Z_{k}$ is a sparse banded matrix. \textit{That is, for a well-conditioned $M$, the dominant entries of the exact solution $Z_{k}$ belong to the sparse banded pattern of the matrix $\underline{Z}_{w}$. Furthermore, as the condition number of $M$ increases, according to the Newton-Schultz error analysis, the parameter $w$ needs to be increased in order to ensure that the \textit{a priori} pattern captures the dominant entries of the exact solution. This implies that the fill-in of the \textit{a priori} pattern $\underline{Z}_{w}$ increases with the condition number of $M$.}

 \textit{The approximation problems for which the matrix $M$ is well-conditioned are referred to as the well-conditioned problems.} Obviously, the condition number of $M$ is closely related to the spectrum or singular values of the matrices $\bar{A}_{k-1}$ and $E$. However, to the best of our knowledge, establishing such a relation is a difficult problem. In the case of non-descriptor state-space models ($E=I$) with symmetric $\bar{A}_{k}=A$, it is easy to show that the condition number of $M$ is equal to the condition number of $A$ \cite{Haber2016sparse}. For small-size systems, we can explicitly form the matrix $M$ and compute its condition number. 
 
\begin{remark}
Although the parameter $w$ is fixed and set by the user, the fill-in of the \textit{a priori} pattern will increase with the increase of the iteration $k$ of the inexact Newton method, summarized in Algorithm~\ref{algorithmNewton}. This is due to the fact that the matrix $\bar{A}_{k-1}$ depends on the previous solution of the GL equation $\hat{Z}_{k-1}$, and in every step $k$ we need to recompute the \textit{a priori} pattern. After a certain number of iterations $k$, the matrix $\bar{A}_{k-1}$ of the inexact Newton method becomes dense, implying a dense \textit{a priori} pattern. To prevent this, the \textit{a priori} pattern should be computed only for small values of $k$, and then such an \textit{a priori} pattern should be fixed and kept constant in further iterations of Algorithm~\ref{algorithmNewton}. 
\label{remarkFillIn}
\end{remark}

\begin{remark}
As stated at the beginning of this section, we are mainly interested in the non-zero structure of the \textit{a priori} pattern, but not in the exact numerical values of its entries. Although the probability of zeros created by coincidence in the numerical values of the matrix entries is very small, it is still a good practice to provide a safeguard against this effect. Furthermore, since in the derivation of the \textit{a priori} pattern,  we ignore the scalar coefficients of the Newton-Schulz iteration, the entries of the matrices $G_{i}$ can rapidly grow or less likely, decrease with $i$, creating numerical instabilities or zeros produced by numerical underflow. Good safeguards against these effects are to operate with the matrices obtained by taking the absolute values (taken element-wise) of the matrices $\bar{A}_{k-1}$, $E$, and $P_{k}$, and after every iteration $i$, to transform the matrices $G_{i}$ into binary matrices according to their non-zero patterns. That is, the $(p,q)$ entry of the binary matrix is equal to $1$ if the corresponding $(p,q)$ entry of $G_{i}$ is non-zero, otherwise, it is zero.
\label{explodingRemark}
\end{remark}

\begin{figure}[H]
\centering
 \includegraphics[scale=0.75,trim=0mm 0mm 0mm 0mm ,clip=true]{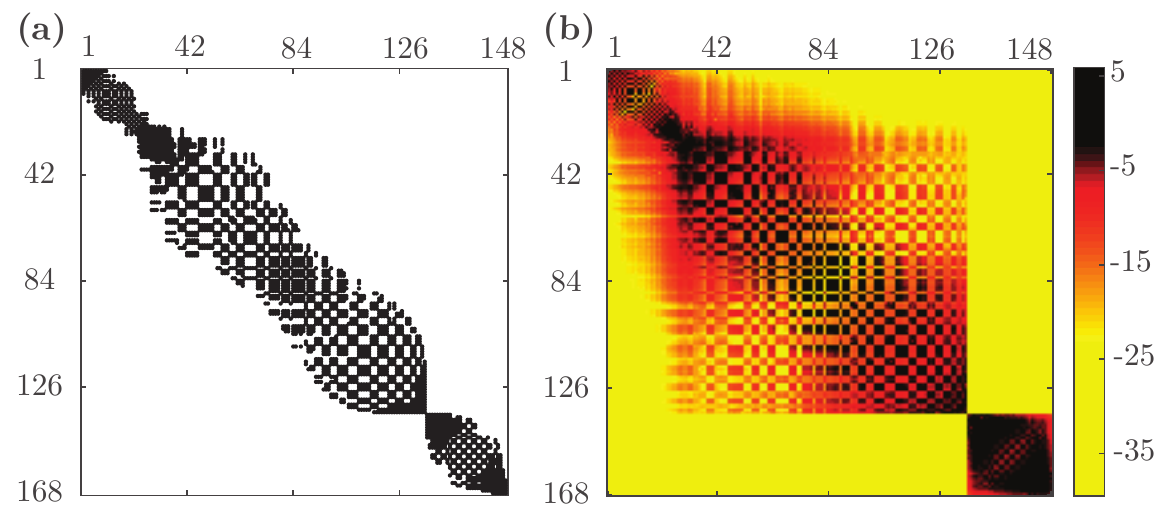}
\caption{The prediction of the pattern of the exact solution of the Lyapunov equation. (a) The pattern predicted by Algorithm~\ref{algorithmApriori}. (b) The absolute values of the entries of the exact solution in the logarithmic scale. The model is obtained using the FE discretization of the heat equation for the domain shown in Fig.~\ref{fig:Gometry}.}
\label{fig:SparsityPrediction}
\end{figure}

Next, we show the numerical performance of the pattern prediction method. Consider the FE model of the heat equation described in Introduction and illustrated in Fig.~\ref{fig:Gometry} (for more details see Section \ref{sectionNumericalExperiments}). We assume that $w=1$ in \eqref{finalAprioriSparsityPattern}, and that the matrices $\bar{A}_{k-1}$ and $P_{k}$ are computed for $\hat{Z}_{k-1}=I$. Figure \ref{fig:SparsityPrediction}(a) shows the predicted pattern of $\hat{Z}_{k}$ whereas Fig.~\ref{fig:SparsityPrediction}(b) shows the absolute values of the entries of the exact solution $Z_{k}$. The exact solution is computed and visually represented using the $\text{lyap}\left(\cdot \right)$ and $\text{imagesc}\left(\cdot \right)$ MATLAB functions.
From Fig.~\ref{fig:SparsityPrediction} we see that the \textit{a priori} pattern is able to accurately predict the most dominant entries of the exact solution.  It is of interest to compute the relevant condition numbers, since according to our theory, they determine the localization degree of the exact solution and the density of the \textit{a priori} pattern. The condition numbers of $E^{T}$ and $\bar{A}_{k-1}^{T}$ are, respectively, $142$ and $18$, whereas the condition number of $M$ is $166$. The matrix $M$ is well-conditioned, implying that the most dominant entries of the exact solution are grouped around the main diagonal, which is numerically confirmed in Fig.~\ref{fig:SparsityPrediction}. 

\section{Methods for Approximating the Solution of the Lyapunov Equation}
\label{sectionLyapunovSolution}
On the basis of the derived \textit{a priori} pattern, in this section we present two methods for approximating the solution of the GL equation. The first method, based on the solution of the reduced normal equations, is a simple generalization of the method presented in \cite{Haber2016sparse}.  For completeness, we briefly summarize this approach in Method~\ref{firstMethod}. However, although the computational complexity of this method is quite low, and for well-conditioned problems it scales almost with $O(n)$, the memory complexity is still significant, regardless of the fact that it still lower than the $O(n^2)$ memory complexity of the MATLAB $\text{lyap}\left(\cdot\right)$ solver. This motivates us to develop a second method based on the Faber expansion of the matrix exponential and on the gradient projection method. This method is a generalization of the Chebyshev expansion method presented in \cite{Haber2016sparse}, which is only applicable to symmetric matrices and NGL equations. The memory complexity of the second method is significantly lower than the first one (for very sparse problems it is $O(n)$), however, the computational complexity is a bit higher, for more details, see Section \ref{sectionNumericalExperiments}. The second approach is summarized in Method~\ref{secondMethod} at the end of this section.

\subsection{First Approximation Method}

Consider the system of equations \eqref{vectorizedLyapunovEquation}. Using the \textit{a priori} pattern of $Z_{k}$ computed using Algorithm~\ref{algorithmApriori}, we can significantly reduce the size of this system and efficiently solve it in the least-squares sense. Let us define $\mathbf{z}_{w}=\text{vec}\left( \underline{Z}_{w}\right)$. The zero entries of the vector $\mathbf{z}_{w}$ determine a zero-entry pattern of the vector that approximates the exact solution $\mathbf{z}$ in \eqref{vectorizedLyapunovEquation}. The dimension of the vector $\mathbf{z}$ should be reduced by eliminating these zero entries. Let the reduced vector be denoted by $\mathbf{z}_{1}\in \mathbb{R}^{n_{1}}$. Consequently, the corresponding columns of the matrix $M$ should be eliminated. More precisely, if the $i$th entry of $\mathbf{z}_{w}$ is zero, then the $i$th column of the matrix $M$ should be eliminated. This column reduction can also produce zero rows of the matrix $M$, and we also assume that such rows are eliminated. Let the reduced matrix be denoted by $M_{1}\in \mathbb{R}^{n_{2}\times n_{1}}$, where $n_{2}\le n^2$ (usually $n_{2}<n^2$), and $n_{1}\ll n^2$, and let the vector $\mathbf{p}$ after possible row reductions be denoted by $\mathbf{p}_{1}\in \mathbb{R}^{n_{2}}$. Once the system of equations has been reduced, the approximate solution can be found by solving the following least-squares problem
\begin{align}
\min_{\mathbf{z}_{1}} \left\| \mathbf{p}_{1}-M_{1}\mathbf{z}_{1} \right\|_{2}^{2}.
\label{leastSquaresProblemsLyapuno}
\end{align}
The solution, denoted by $\hat{\mathbf{z}}_{1}$, can be found by solving the system of normal equations
\begin{align}
M^{T}_{1}M_{1}\hat{\mathbf{z}}_{1}=M_{1}^{T}\mathbf{p}_{1},
\label{normalEquations}
\end{align}
using for example the Conjugate Gradient Least Squares (CGLS) method~\cite{bjorck1996numerical}. The computational and memory complexities of the first method are thoroughly analyzed and numerically simulated in \cite{Haber2016sparse}, and for brevity we omit this analysis. In brief, one iteration of the CGLS costs approximately $\text{nnz}\left(M_{1}\right)+3n_{1} +2n_{2}$ flops, where $\text{nnz}\left(\cdot\right)$ denotes the number of non-zero entries determined by the parameter $w$. If $M$ is well-conditioned, $w$ can be chosen as a small number, and the solution can be computed with (almost) $O(n)$ computational complexity. \textit{However, special attention needs to be given to the implementation of the first method.} Namely, a naive implementation based on a direct formation of the matrix $M$ using for example, the MATLAB $\text{kron}\left(\cdot\right)$ function, followed by a reduction of this matrix, is not possible for a large $n$ due to the high memory demand. Instead, the non-zero entries of the reduced matrix $M_{1}$ should be calculated without forming the matrix $M$. This itself requires some computational time, however, our numerical experience shows that with proper implementation in the MATLAB sparse matrix toolbox, this time is shorter than the time necessary to solve the least-squares problem. Also, although for well-conditioned problems the memory complexity scales almost linearly with $n$, the reduced matrix $M_{1}$ for a large $n$ still occupies significant memory space, see Section \ref{sectionNumericalExperiments} for more details.

\setcounter{algorithm}{0}
\begin{method}[H]
\caption{Approximation of the GL equation using the reduced least-squares problem}
\begin{enumerate}
\item Compute the \textit{a priori} pattern $\underline{Z}_{w}$ using Algorithm~\ref{algorithmApriori}.
\item On the basis of the non-zero structure of $\underline{Z}_{w}$, form the reduced matrix $M_{1}$ and solve the least-squares problem \eqref{leastSquaresProblemsLyapuno} using the CGLS.
\end{enumerate}
\label{firstMethod}
\end{method}

\subsection{Second Approximation Method}

The approximate solution of the GL equation \eqref{lyapunovEquationReWritten} can be found by solving the following optimization problem:
\begin{align}
& \min_{Z_{k}}\left\|P_{k}-E^{T}Z_{k}\bar{A}_{k-1}-\bar{A}_{k-1}^{T}Z_{k}E \right\|_{F}^{2}, \;\; \text{subject to}\;\; Z_{k} \in \Omega_{\underline{Z}_{w}},
\label{solutionLyapunovOptimization} 
\end{align}
where $\Omega_{\underline{Z}_{w}}$ is a set of all matrices with the pattern defined by the pattern of the matrix $\underline{Z}_{w}$ computed by Algorithm~\ref{algorithmApriori}. The solution of \eqref{solutionLyapunovOptimization} can be obtained using the \textit{gradient projection method} \cite{bertsekas1999nonlinear}:
\begin{align}
\tilde{Z}_{k}^{i+1}=\mathcal{V}_{\underline{Z}_{w}}\big[\tilde{Z}_{k}^{i}-\delta^{i}N_{k}^{i} \big],
\label{gradientProjection}
\end{align}
where $\tilde{Z}_{k}^{i}$ is the approximate solution at the $i$-th iteration, $\mathcal{V}_{\underline{Z}_{w}}\left(\cdot\right)$ is the projection operator defined in \eqref{sparsificationOperatordefinition} (this operator projects the argument matrix on the pattern of $\underline{Z}_{w}$), $N_{k}^{i}$ is the gradient matrix, and $\delta^{i}$ is the step size. Expressing the cost function \eqref{solutionLyapunovOptimization} as a sum of matrix traces, and by differentiating such an expression, it can be shown that
\begin{align}
& N_{k}^{i}=-2ER_{k}^{i}\bar{A}_{k-1}^{T}-2\bar{A}_{k-1}R_{k}^{i}E^{T}, \;\; R_{k}^{i}=P_{k}-E^{T}\tilde{Z}_{k}^{i}\bar{A}_{k-1}-\bar{A}_{k-1}^{T}\tilde{Z}_{k}^{i}E.
\label{gradientExplained}
\end{align}
To define the step size $\delta^{i}$, we first introduce
\begin{align}
& J\big[\tilde{Z}_{k}^{i}\big]=\left\|P_{k}-E^{T}\tilde{Z}_{k}^{i}\bar{A}_{k-1}-\bar{A}_{k-1}^{T}\tilde{Z}_{k}^{i}E \right\|_{F}^{2}, \; \tilde{Z}_{k}^{i}\big[\delta \big]=\mathcal{V}_{\underline{Z}_{w}}\big[\tilde{Z}^{i}_{k}-\delta N_{k}^{i} \big]. \label{stepSize1}
\end{align}
Then, the step size is determined by the Armijo rule along the projection arc, defined by $\delta^{i}=\zeta^{g_{i}}\overline{\delta}$, where $\overline{\delta}>0$, $\zeta\in (0,1)$, and $g_{i}$ equals the first nonnegative integer $g$ satisfying
\begin{align}
J\big[\tilde{Z}_{k}^{i}\big]-J\Big[\tilde{Z}_{k}^{i}\big[ \zeta^{g}\overline{\delta}\big]\Big]\ge \sigma\text{vec}\left(N_{k}^{i}\right)^{T}\text{vec}\left(\tilde{Z}_{k}^{i}-\tilde{Z}_{k}^{i}\big[\zeta^{g}\overline{\delta} \big] \right),
\label{armijoRule}
\end{align}
and where $\sigma>0$. The convergence of the gradient projection method can be analyzed by recognizing that the cost function in \eqref{solutionLyapunovOptimization} can be written in the vectorized form, which corresponds to the least squares formulation of the solution to the linear system of equations \eqref{vectorizedLyapunovEquation}. From there it follows that the convergence rate is determined by the minimal and maximal singular values of $M$, that is, it is determined by the condition number of $M$ \cite{bertsekas1999nonlinear}. Provided that the parameter $w$ is much smaller than $n$, one step of the gradient projection method can be implemented with $O(n)$ computational and memory complexities. In the case of ill-conditioned problems, the convergence of the gradient projection method can be slow. \textit{Consequently, for fast convergence, it is essential that it be initialized with a good initial guess of the solution.} In the sequel we present a method for generating such an initial guess. By first multiplying \eqref{lyapunovEquationReWritten} from the left with $\left(E^{T}\right)^{-1}$ and then from the right with $E^{-1}$, we transform the GL equation into the NGL equation:
\begin{align}
Z_{k}\bar{A}_{k-1}E^{-1}+\left(E^{T}\right)^{-1}\bar{A}_{k-1}^{T}Z_{k}=\left(E^{T}\right)^{-1}P_{k}E^{-1}.
\label{multiplicationLeftRightLyapunov}
\end{align}
Assuming that $\bar{A}_{k-1}E^{-1}$ is a stable matrix, the unique solution of \eqref{multiplicationLeftRightLyapunov} can be expressed as follows (see Theorem 13.23 in \cite{laub2005matrix}):
\begin{align}
\mathcal{X}=-\int_{0}^{\infty}\exp\left(t \mathcal{A}\right) \mathcal{P} \exp\left(t \mathcal{A}^{T}\right) \mathrm{d}t,
\label{integralRepresentationLyapunov}
\end{align}
where $\mathcal{A}=\left(E^{-1} \right)^{T} \bar{A}^{T}_{k-1}$, $\mathcal{P}=\left(E^{-1}\right)^{T}P_{k}E^{-1}$, and where we rely upon the fact that $\left( E^{-1}\right)^{T}=\left( E^{T}\right)^{-1}$. By inverting the matrix $E$, the sparsity of the problem has been lost. To preserve the sparsity, we search for a sparse approximate inverse of $E$. Such an inverse exists due to the following facts. Given that the matrix $E$ is banded, under some mild conditions its inverse will exhibit some form of localization, or in the positive-definite case it will belong to the class of off-diagonally decaying matrices, see for example~\cite{demko1984,benzi2007,benzi2017approximation,benzi2015decay,benzi2014decay}. In the general case, when the matrix $E$ is not strictly banded or positive-definite, the sparse approximation of $E^{-1}$ still exists if the matrix $E$ is relatively well-conditioned. This follows from the convergence analysis of the Newton-Schulz method, see Section~\ref{sectionAprioriSparsityPattern}. All these facts motivate us to search for a sparse banded approximate inverse of $E$. This approximate inverse can be computed by solving the following optimization problem \cite{grote1997parallel,chow1998approximate}
\begin{align}
\min_{\mathcal{E}}\left\|I-E\mathcal{E} \right\|_{F}^{2}, \;\; \text{subject to}\;\; \mathcal{E}\in \Omega_{\mathcal{E}}.
\label{optimizationM}
\end{align} 
where $\mathcal{E}$ is a sparse approximate inverse of $E$ and $\Omega_{\mathcal{E}}$ is a set of matrices with a given pattern. The problem \eqref{optimizationM} is highly parallelizable and can be split into $n$ independent low-dimensional least-squares problems that can be efficiently solved using the SPAI solver~\cite{grote1997parallel,chow1998approximate,tang1999toward,chow2000priori,benzi1998sparse}. In some situations, it might be more appropriate to replace the term $I-E\mathcal{E}$ by $I-\mathcal{E}E$ in \eqref{optimizationM}. In practice, both forms can be easily computed and the most accurate solution should be selected.
The choice of the set  $\Omega_{\mathcal{E}}$ depends on the pattern of $E$. Using arguments analogous to those in Section \ref{sectionAprioriSparsityPattern}, it can be shown that a good estimate of the pattern of $\mathcal{E}$ is given by $I+E+E^2+\ldots+E^{k_{1}}$, where $k_{1}$ is a small positive integer.  Similarly to the analysis developed in Section \ref{sectionAprioriSparsityPattern}, it can be shown that $k_{1}$ is generally small for a well-conditioned $E$. By substituting in \eqref{integralRepresentationLyapunov}, the matrix $E^{-1}$ by the matrix $\mathcal{E}$, that solves \eqref{optimizationM}, we define the following approximation 
\begin{align}
\mathcal{X} \approx \mathcal{X}_{1}, \;\; \mathcal{X}_{1}=-\int_{0}^{\infty}\exp\left(t \mathcal{A}_{1}\right) \mathcal{P}_{1}\exp\left(t \mathcal{A}_{1}^{T}\right) \mathrm{d}t, \;\; 
\label{integralRepresentationLyapunovApproximation}
\end{align}
where $\mathcal{A}_{1}=\mathcal{E}^{T}\bar{A}_{k-1}^{T}$ and $\mathcal{P}_{1}=\mathcal{E}^{T}P_{k}\mathcal{E}$. Because the matrices $\mathcal{E}$, $\bar{A}_{k-1}$, and $P_{k}$ are banded, the matrices $\mathcal{A}_{1}$ and $\mathcal{P}_{1}$  are also banded. 

Our goal is to approximate $\mathcal{X}_{1}$ by a banded matrix. Denote with $\lambda_{\text{RS}}$, $\lambda_{\text{RL}}$, and $\lambda_{\text{IL}}$ the smallest real part, the largest real part, and the largest absolute value of the imaginary part of eigenvalues of $\mathcal{A}_{1}$. Because the matrix $\mathcal{A}_{1}$ is sparse banded, its extreme eigenvalues can be estimated with $O(n)$ complexity using the ARPACK software (for more details see Section 4 of \cite{bergamaschi2003efficient}). From Theorem 4.3 in \cite{grasedyck2003solution}, we have that the integral \eqref{integralRepresentationLyapunovApproximation} can be approximated by:
\begin{align}
\mathcal{X}_{2}=-\sum_{j=-q}^{q} \psi \omega_{j}\exp\left( \tilde{t}_{j} \mathcal{A}_{1} \right)\mathcal{P}_{1}\exp\left(  \tilde{t}_{j} \mathcal{A}_{1}^{T} \right),
\label{method1LyapunovApproximate}
\end{align}
where $\tilde{t}_{j}=\psi t_{j}$, $q$ is a sufficiently large integer, and 
\begin{align}
&\psi=\frac{3}{2|\lambda_{RL}|}, \; \; \omega_{j}=\left(q+q\exp\left(-2jq^{-1/2} \right)\right)^{-1/2}, \notag \\
&t_{j}=\log\left( \exp\left(jq^{-1/2} \right)+\sqrt{1+\exp\left(2jq^{-1/2}\right)} \right). 
\label{method1LyapunovApproximateExplanation}
\end{align}
In \cite{grasedyck2003solution} it has been shown that the approximation error exponentially decreases with $\sqrt{q}$, that is, $\left\|\mathcal{X}_{1}-\mathcal{X}_{2}\right\|_{2}\le K_{\mathcal{A}_{1}}\left\|\mathcal{P}_{1}\right\|_{2}\exp\left(-\sqrt{q} \right)$, where the constant $K_{\mathcal{A}_{1}}$ depends on $\lambda_{RS}$, $\lambda_{RL}$ and $\lambda_{IL}$. 

To compute $\mathcal{X}_{2}$, we need to approximate the matrix exponential by a banded  matrix. Consequently, in the sequel we develop an algorithm for approximating the matrix exponential by a banded matrix. The developed approach is summarized in Algorithm~\ref{algorithmExponential}. First of all, the justification of the existence of the banded approximation of the matrix exponential follows from the following facts. If the matrix $\tilde{t}_{j}\mathcal{A}_{1}$ is diagonalizable for each $\tilde{t}_{j}\in \left[0,\infty \right]$, then its matrix exponential is an off-diagonally decaying matrix \cite{benzi2007}. In \cite{pozza2016decay} these results are extended to non-Hermitian, not necessarily diagonalizable matrices. Furthermore,  our numerical results indicate that the matrix exponentials of a wide class of sparse banded matrices belong to the class of localized matrices.  

We use the Faber polynomials to approximate the matrix exponential \cite{bergamaschi2003efficient}. Denote with $\tilde{\lambda}_{\text{RS},j}$, $\tilde{\lambda}_{\text{RL},j}$, and $\tilde{\lambda}_{\text{IL},j}$ the smallest real part, the largest real part, and the largest  imaginary part of the eigenvalues of  $\tilde{t}_{j}\mathcal{A}_{1}$ (these scalars can be computed directly from $\lambda_{\text{RS}}$, $\lambda_{\text{RL}}$, and $\lambda_{\text{IL}}$). Next, define \cite{bergamaschi2003efficient}:
\begin{align}
&c_{1,j}=\frac{\tilde{\lambda}_{\text{RL},j}-\tilde{\lambda}_{\text{RS},j}}{2}, \; c_{2,j}=\frac{c_{1,j}^{\frac{2}{3}}\sqrt{c_{1,j}^{\frac{2}{3}}+\tilde{\lambda}_{\text{IL},j}^{\frac{2}{3}}}+\sqrt{\left(c_{1,j}\tilde{\lambda}_{\text{IL},j}^{2} \right)^{\frac{2}{3}}+\tilde{\lambda}_{\text{IL},j}^{2}}}{2}, \notag \\
&c_{3,j}=\left(c_{1,j}^{\frac{2}{3}}+\tilde{\lambda}_{\text{IL},j}^{\frac{2}{3}}\right)\left(c_{1,j}^{\frac{4}{3}}-\tilde{\lambda}_{\text{IL},j}^{\frac{4}{3}} \right), \; c_{4,j}=\frac{\tilde{\lambda}_{\text{RL},j}+\tilde{\lambda}_{\text{RS},j}}{2}.
\label{definitionEigenvaluesConstants}
\end{align}
The matrix $\mathcal{A}_{2,j}$ is defined as follows \cite{bergamaschi2003efficient}:
\begin{align}
\mathcal{A}_{2,j}=\frac{1}{\sqrt{c_{3,j}}}\left(\tilde{t}_{j}\mathcal{A}_{1}-c_{4,j}I\right).
\label{matrixA3definition}
\end{align}
The truncated Faber series approximation of  $\exp \left(\tilde{t}_{j} \mathcal{A}_{1}\right)$ is defined by \cite{bergamaschi2003efficient}:
\begin{align}
\mathcal{K}_{p,j}=\sum_{l=0}^{p}a_{l,j}F_{l,j},
\label{faberSeries}
\end{align}
where $p$ is a positive integer, Faber coefficients $a_{l,j}$ are defined in Appendix, and the Faber polynomials $F_{l,j}$ are defined by:
\begin{align}
& F_{0,j}=I, \; F_{l,j}=2\left(\frac{\sqrt{c_{3,j}}}{2c_{2,j}} \right)^{l}T_{l,j}\left( \mathcal{A}_{2,j} \right), \;\; l>0, \label{faberPoly}
\end{align}
and where $T_{l,j}\left( \mathcal{A}_{2,j} \right)$ are the Chebyshev polynomials of the first kind defined by the recurrence relation:
\begin{align}
&T_{0,j}=I, \; T_{1,j}=\mathcal{A}_{2,j}, \; T_{k+1,j}=2\mathcal{A}_{2,j}T_{k,j}-T_{k-1,j},\; k=1,2,\ldots
\label{chebyshevPoly}
\end{align}
An upper bound on the approximation error of the truncated Faber series is given in Section 3 of \cite{bergamaschi2003efficient} and for brevity we omit it. Briefly speaking, the approximation error is smaller if $p$ is lager, however, larger values of $p$ decrease the number of non-zero elements of $\mathcal{K}_{p,j}$. This is due to the fact that the recurrence \eqref{chebyshevPoly} implies that the matrix $\mathcal{K}_{p,j}$ in \eqref{faberSeries} can be expressed as the sum of the first $p$ powers of $\mathcal{A}_{2,j}$. If $p$ is much smaller than $n$, then $\mathcal{K}_{p}$ is a sparse banded matrix and it can be computed with $O(n)$ complexity \cite{benzi2007}. However, if this condition is not satisfied, then the matrix $\mathcal{K}_{p,j}$ is no longer sparse and consequently, it cannot be computed with $O(n)$ complexity. To ensure that  $\mathcal{K}_{p,j}$ remains sparse for larger values of $p$, and to ensure that it can be computed with $O(n)$ complexity, we sparsify the recurrence relation \eqref{chebyshevPoly} as follows:
\begin{align}
&\tilde{T}_{0,j}=I, \tilde{T}_{1,j}=\mathcal{A}_{2,j}, \; \tilde{T}_{k+1,j}=\mathcal{V}_{\underline{\mathcal{A}}_{2}}\big[2\mathcal{A}_{2,j}\tilde{T}_{k,j}-\tilde{T}_{k-1,j}\big],\; k=1,2,\ldots
\label{chebyshevPolySparsified}
\end{align}
where $\mathcal{V}_{\underline{\mathcal{A}}_{2}}\big[\cdot\big]$ is a projection operator constructed for example by summing the identity matrix with the first few powers of $\mathcal{A}_{2,j}$ (or $\mathcal{A}_{1}$). In this way, we also sparsify the Faber polynomials. 
Let the sparsified approximation of the matrix exponential, obtained by using the sparsified Chebyshev polynomials \eqref{chebyshevPolySparsified}, be denoted by $\tilde{\mathcal{K}}_{p,j}$
\begin{align}
\tilde{\mathcal{K}}_{p,j}=\sum_{l=0}^{p}a_{l,j}\tilde{F}_{l,j},
\label{faberSeriesSparsfied}
\end{align}
where $\tilde{F}_{l,j}$ are the sparsified Faber polynomials computed by substituting \eqref{chebyshevPolySparsified} in \eqref{faberPoly}.

\setcounter{algorithm}{2}
\begin{algorithm}[H]
\caption{Computation of the sparse banded approximation of the matrix exponential $\exp\left(\tilde{t}_{j} \mathcal{A}_{1} \right)$ in \eqref{method1LyapunovApproximate}}
\textbf{Input:} the matrices $E$, $\bar{A}_{k-1}$, parameter $\tilde{t}_{j}$, and the approximation order $p$ in \eqref{faberSeries}. \\
\textbf{Output:} the approximation of the matrix exponential $\tilde{\mathcal{K}}_{p,j}$.
\begin{enumerate}
\item  Computation of the sparse approximate inverse $\mathcal{E}$ of the matrix $E$. 
\begin{enumerate}
\item Compute the matrix $I+E+E^2+\ldots+E^{k_{1}}$ that represents an \textit{a priori} pattern of $E^{-1}$. The small positive integer $k_{1}$ is a user choice, and it should be selected such that the a priory pattern is sparse ($k_{1}\ll n$).
\item For the computed \textit{a priori} pattern, compute the approximate inverse $\mathcal{E}$ by solving  \eqref{optimizationM}  using the SPAI solver. 
\end{enumerate}
\item Compute the matrix  $\mathcal{A}_{1}=\mathcal{E}^{T}\bar{A}_{k-1}^{T}$.
\item Compute $\tilde{\lambda}_{\text{RS},j}$, $\tilde{\lambda}_{\text{RL},j}$, and $\tilde{\lambda}_{\text{IL},j}$ that represent the smallest real part, the largest real part, and the largest imaginary part of the eigenvalues of $\tilde{t}_{j}\mathcal{A}_{1}$. These scalars can be computed from the extreme parts of the eigenvalues of $\mathcal{A}_{1}$.
\item Compute the constants $c_{1,j}$, $c_{2,j}$, $c_{3,j}$ and $c_{4,j}$ defined in \eqref{definitionEigenvaluesConstants} and compute the matrix $\mathcal{A}_{2,j}$ defined in \eqref{matrixA3definition}.
\item Compute the coefficients of the truncated Faber series expansion (the coefficients are defined in Appendix).
\item  Computation of the sparsified Faber polynomials $\tilde{F}_{l,j}$.
\begin{enumerate}
\item  Compute the pattern of the projection operator  $\mathcal{V}_{\underline{\mathcal{A}}_{2}}\big[\cdot\big]$ of the sparsified Chebyshev polynomials in \eqref{chebyshevPolySparsified}. The pattern is $I+\mathcal{A}_{2,j}+\ldots \mathcal{A}_{2,j}^{k_{2}}$, where $k_{2}$ is a small positive integer ($k_{2}\ll n$) that should be chosen such that the resulting matrix remains sparse.
\item Compute the sparsified Chebyshev polynomials \eqref{chebyshevPolySparsified}. Substitute such polynomials in \eqref{faberPoly} and compute the sparsified Faber polynomials.
\end{enumerate}
\item Compute the approximation of the matrix exponential \eqref{faberSeriesSparsfied}.
\end{enumerate}
\label{algorithmExponential}
\end{algorithm}

 Substituting $\tilde{\mathcal{K}}_{p,j}$ in \eqref{method1LyapunovApproximate}, we obtain:
\begin{align}
& \mathcal{X}_{3}=-\sum_{j=-q}^{q} \psi \omega_{j}  \tilde{\mathcal{K}}_{p,j}\mathcal{P}_{1}\Big(\tilde{\mathcal{K}}_{p,j}\Big)^{T}. \label{finalApproximationFaber}
\end{align}
The matrix \eqref{finalApproximationFaber} is used as an initial guess of the gradient projection method \eqref{gradientProjection}. The initial guess \eqref{finalApproximationFaber} can be computed with $O(n)$ computational and memory complexities provided that $q\ll n$, $p$ in \eqref{faberSeries} satisfies $p\ll n$, and that the projection operator in \eqref{chebyshevPolySparsified} is chosen such that the projected matrix remains sparse. In practice, we do not need a high accuracy of the initial guess, and consequently, some inaccuracies can be tolerated for the benefit of computing it with low computational complexity. The second approximation method is summarized in Method~\ref{secondMethod}.

\setcounter{algorithm}{1}
\begin{method}[H]
\caption{Approximation of the GL equation using the gradient projection method and the sparsified Chebyshev polynomials}
\begin{enumerate}
\item Compute the \textit{a priori} pattern $\underline{Z}_{w}$ using Algorithm~\ref{algorithmApriori}.
\item Generate an initial guess for the gradient projection method.
\begin{enumerate}
\item Choose a relatively small integer $q$ in \eqref{finalApproximationFaber} ($q\ll n $). Compute $\mathcal{P}_{1}=\mathcal{E}^{T}P_{k}\mathcal{E}$. Compute the parameters $\tilde{t}_{j}$, $\psi$, and $\omega_{j}$ using \eqref{method1LyapunovApproximateExplanation}.
\item For $j=-q,-q+1,\ldots, q$, compute the matrix $\tilde{\mathcal{K}}_{p,j}$ using Algorithm~\ref{algorithmExponential}. Compute the initial guess \eqref{finalApproximationFaber}.
\end{enumerate}
\item Set  $\tilde{Z}_{k}^{0}=\mathcal{X}_{3}$, and propagate the gradient projection iteration \eqref{gradientProjection} until convergence or until the predefined maximal number of iterations has been reached. 
\end{enumerate}
\label{secondMethod}
\end{method}

\begin{remark}
We have developed the approximation methods assuming sparse banded structure of the system matrices $E$ and $A$. Here we briefly explain the generalization of the proposed methods to the case of model \eqref{descriptorStateSpaceModel} where the matrices $E$ and $A$ are composed of sparse banded blocks. The simplest choice of the initial guess in Algorithm~\ref{algorithmNewton} is a $2 \times 2 $ block matrix (blocks corresponding to the blocks of the  matrices $E$ and $A$ in \eqref{descriptorStateSpaceModel}), in which every block is an identity matrix.  Assuming such an initial guess, the \textit{a priori} pattern computed using Algorithm \ref{algorithmApriori}, is a $2\times 2$ block matrix with (sparse) banded blocks. This is because the pattern is computed by multiplying the matrices $E$, $A$, and $P_{k}$, whose blocks are (sparse) banded matrices. The fill-in of these blocks will be increased after every multiplication, however, as before, we assume that only a relatively small number of multiplications is necessary, and consequently, the blocks will remain sparse. Once this pattern has been computed, Method~\ref{firstMethod} and Method~\ref{secondMethod} can be directly applied. 
\label{remarkGenerati}
\end{remark}

\section{Numerical Experiments}
\label{sectionNumericalExperiments}
\begin{figure}[H]
\centering
 \includegraphics[scale=0.4,trim=0mm 0mm 0mm 0mm ,clip=true]{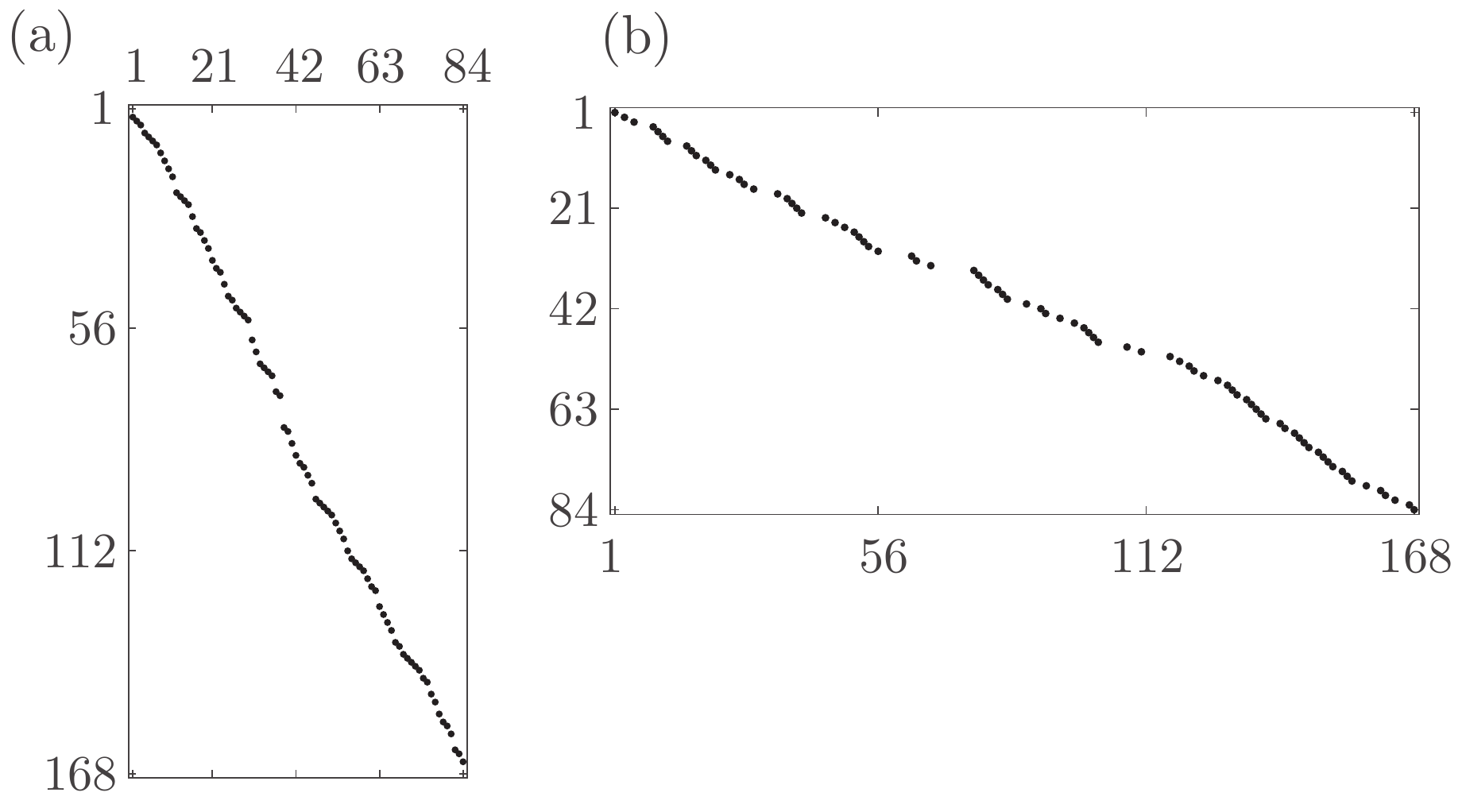}
\caption{The patterns of (a) the matrix $B$ and (b) the matrix $C$.}
\label{fig:SparsityBC}
\end{figure}

\begin{figure}[H]
\centering
 \includegraphics[scale=0.51,trim=0mm 0mm 0mm 0mm ,clip=true]{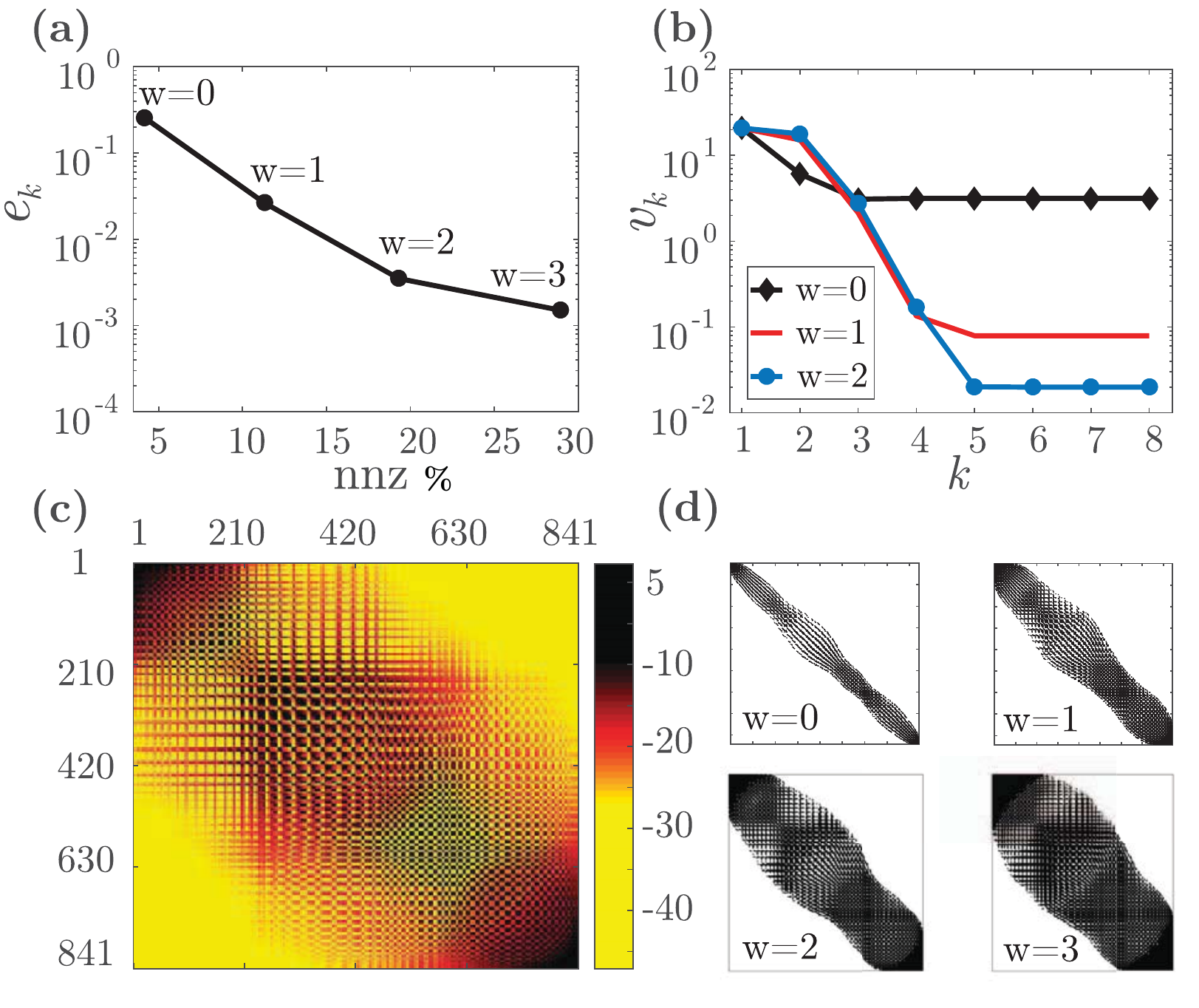}
\caption{(a) The accuracy of solving the GL equation, measured by $e_{k}$ defined in \eqref{eLyapunov}, using the first method for several \textit{a priori} patterns. (b) The convergence of the inexact Newton method, measured by $v_{k}$ defined in \eqref{residualRiccati}, for $w=0,1,2$. (c) An image of the exact solution of the GR equation. The colors of the pixels correspond to the logarithms of the absolute values of matrix entries. (d)\textit{A priori} patterns for several values of $w$.}
\label{fig:ConvergenceSparsity}
\end{figure}

\begin{figure}[H]
\centering
 \includegraphics[scale=0.25,trim=0mm 0mm 0mm 0mm ,clip=true]{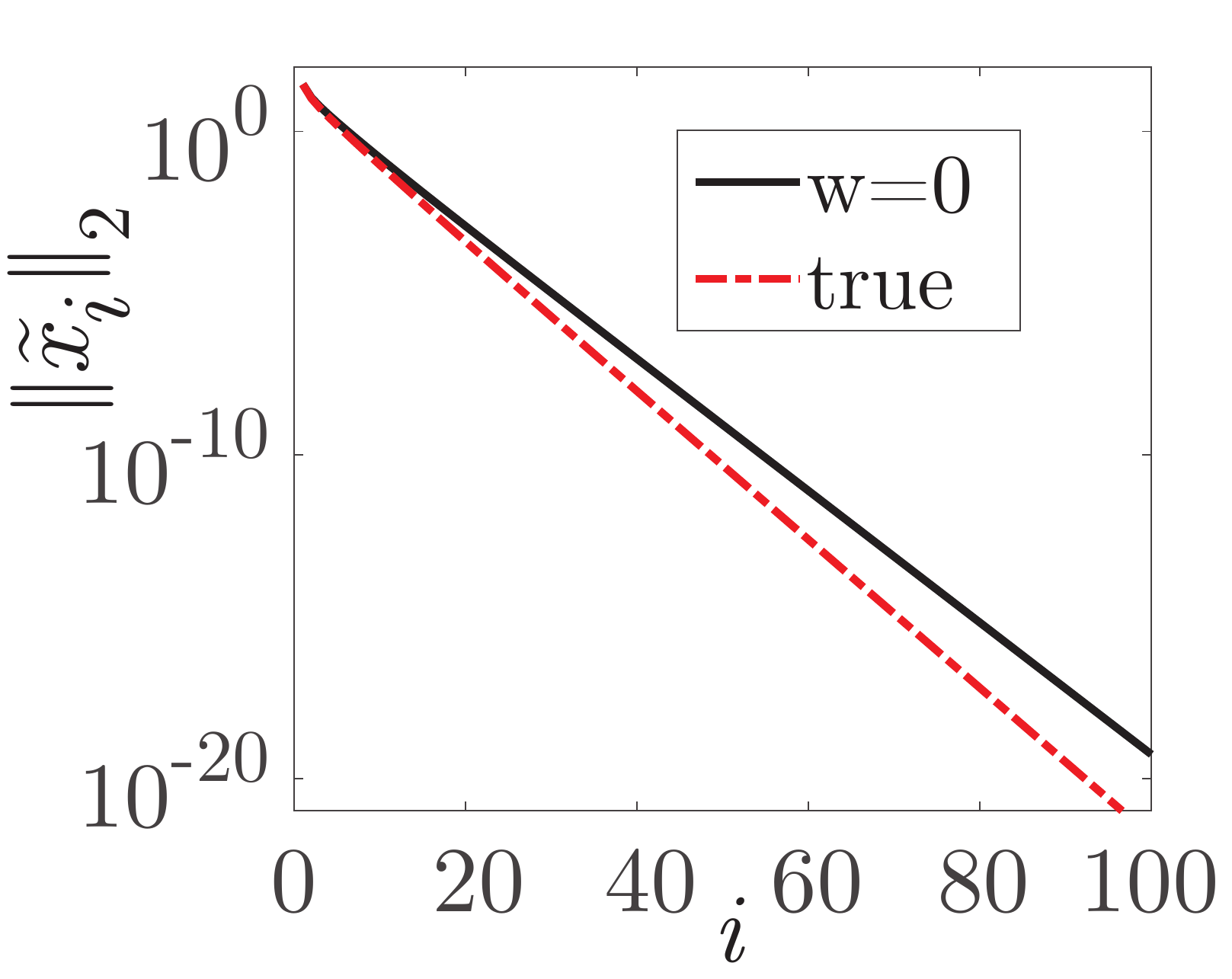}
\caption{The closed loop performance of the system for the exact feedback matrix (denoted by "true") and for the approximate feedback matrix computed using the first method for $w=0$. The integer $i$ is the discrete-time index used in the simulation of the continuous time dynamics.}
\label{fig:PerformanceControl}
\end{figure}

The simulations are performed on a desktop computer with $16$GB RAM and Intel Xeon Processor E3-1245 v5. We use the FE model of the discretized heat equation described in Introduction. The MATLAB codes and the models used in our simulations can be found in~\cite{haber2018codes}. The discretization domain is shown in Fig.~\ref{fig:Gometry}(a). We generate $4$ models with an increasing number of nodes. As the node number increases, the mesh becomes denser, and the state dimension increases. We use the models having $n=168,841,3687$ and $15379$ states. The matrix $B\in \mathbb{R}^{n\times \lfloor n/2 \rfloor}$ is constructed by randomly placing the localized actuators on $\lfloor n/2 \rfloor$ nodes in the 2D domain in Fig.~\ref{fig:Gometry}(a). Similarly, we construct the matrix $C\in \mathbb{R}^{\lfloor n/2 \rfloor \times n}$. The patterns of the matrices $B$ and $C$ are shown in Fig.~\ref{fig:SparsityBC}. For simplicity, the weighting matrices $Q$ and $R$ in \eqref{theCostFunction} are chosen as $Q=I$ and $R=I$, where $I$ is an identity matrix. The initial guess for the solution of the Riccati equation in Algorithm~\ref{algorithmNewton} is $\hat{Z}_{0}=10 I$. This is a reasonable initial guess, because it is expected that the exact solution is a localized matrix with the localization region around the main diagonal. 
 The exact solutions of the Lyapunov and Riccati equations are computed using the built-in MATLAB functions $\text{lyap}\left(\cdot \right)$ and $\text{dare}\left(\cdot \right)$, respectively. Due to the $O(n^3)$ computational and $O(n^2)$ memory complexities of the algorithms implemented in these functions, we are only able to compute the exact solutions for the models for which $n\le 3687$. The accuracy of the proposed methods for solving the Lyapunov equations is quantified by comparing it to the exact solution
\begin{align}
e_{k}=\left\|\hat{Z}_{k}-Z_{k}\right\|_{F}/\left\|Z_{k}\right\|_{F}.
\label{eLyapunov}
\end{align} 
 The convergence of the inexact Newton method is quantified using the Frobenius norm of the residual
\begin{align}
v_{k}=\left\| \mathcal{D}\left[\hat{Z}_{k} \right]\right\|_{F},
\label{residualRiccati}
\end{align}
where the Riccati residual $\mathcal{D}\left[\cdot \right]$ is defined in \eqref{riccatiOperator}. First, we show the accuracy of the approximate solution of the GL equation as a function of the number of nonzero entries of the \textit{a priori} pattern $\underline{Z}_{w}$. The results of solving the GL equation in the first iteration of the inexact Newton method (Algorithm~\ref{algorithmNewton}) are shown in Fig.~\ref{fig:ConvergenceSparsity}(a), for $w=0,1,2,3$. The GL equation is solved using the first method. The residual convergence tolerance of the CGLS method is $10^{-7}$. The \textit{a priori} patterns are shown in Fig.~\ref{fig:ConvergenceSparsity}(d). Figure~\ref{fig:ConvergenceSparsity}(b) shows the convergence of the inexact Newton method for solving the GR equation for $w=0,1,2$. It can be observed that the convergence is relatively stable and that the steady-state residual decreases as the number of the nonzero entries of the \textit{a priori} pattern increases. This is a direct consequence of the fact that the accuracy of solving the GL equation, shown in Fig.~\ref{fig:ConvergenceSparsity}(a), increases with an increase of the number of nonzero entries of the \textit{a priori} patterns. Figure~\ref{fig:ConvergenceSparsity}(c) shows the exact solution of the GR equation. This figure confirms that the exact solution belongs to the class of localized matrices and confirms that the \textit{a priori} patterns shown in Fig.~\ref{fig:ConvergenceSparsity}(d) are able to accurately capture the most dominant entries of the exact solution. 
\begin{figure}[H]
\centering
 \includegraphics[scale=0.48,trim=0mm 0mm 0mm 0mm ,clip=true]{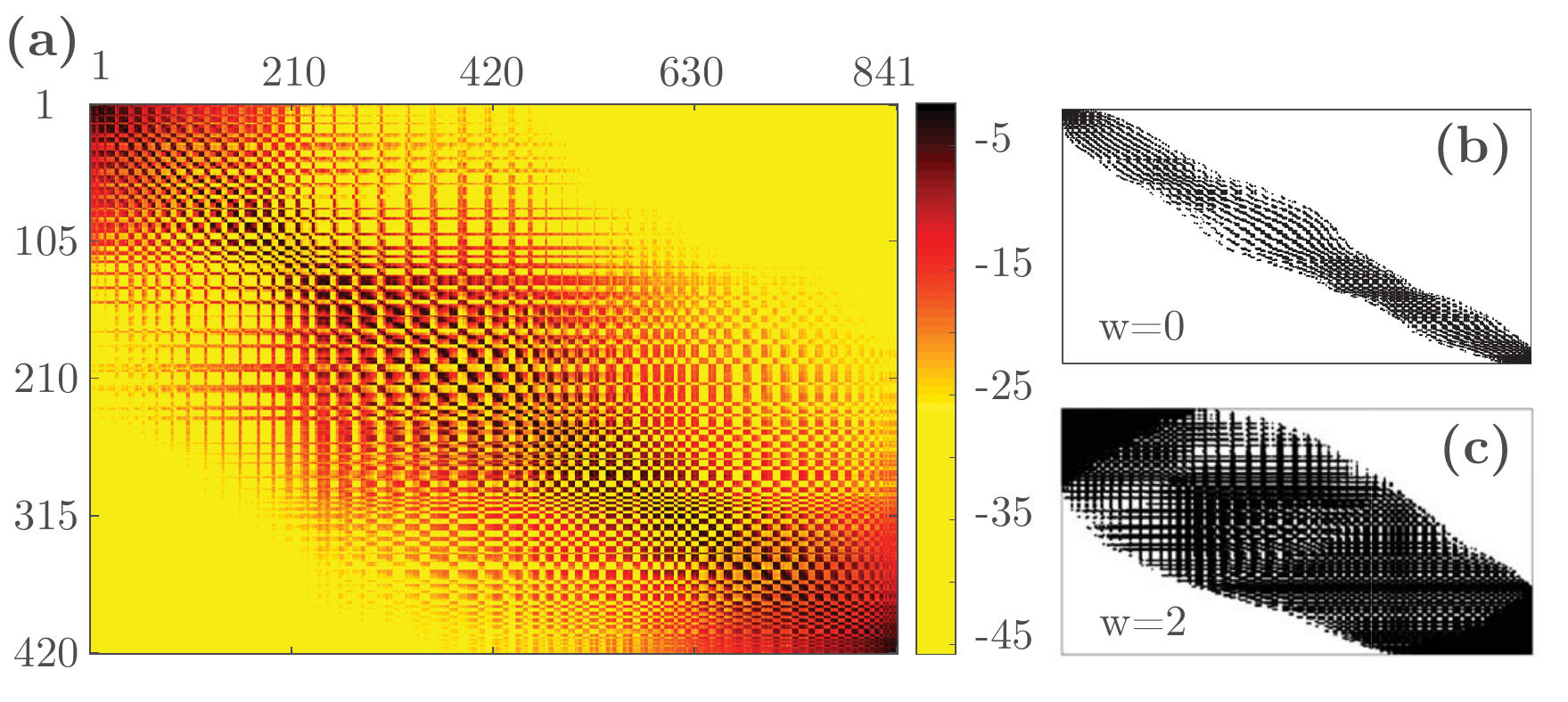}
\caption{(a) The absolute values of the entries of the exact feedback matrix $F$ in the logarithmic scale. The patterns of the computed feedback matrix for the \textit{a priori} patterns determined by (b) $w=0$ and (c) w=2.}
\label{fig:FeedbackMatrixF}
\end{figure}
Figure~\ref{fig:FeedbackMatrixF}(a) shows the feedback control matrix $F$ computed on the basis of the exact solution of the GR equation. Figures \ref{fig:FeedbackMatrixF}(b) and (c) show the patterns of the approximate feedback control matrices $\hat{F}$ computed for $w=0$ and $w=2$, respectively. The percentages of non-zero entries of $\hat{F}$, for $w=0,1,2$ and $3$ are, respectively, $5.7,13.0,21.4$ and $31.3$. It can be observed that the feedback control matrices computed for $w=0$ and $w=2$ are able to accurately capture the most dominant entries of the exact feedback matrix. We compare the performance of the closed loop system controlled using the exact feedback matrix and the approximate feedback matrix calculated using the first method. The results are shown in Fig.~\ref{fig:PerformanceControl}. It can be observed that the loss of performance of the approximate control law is almost negligible compared to the exact control law. Notice that the presented results are generated for $w=0$, implying that for $w>0$, performance will be even better.

\begin{figure}[H]
\centering
 \includegraphics[scale=0.4,trim=0mm 0mm 0mm 0mm ,clip=true]{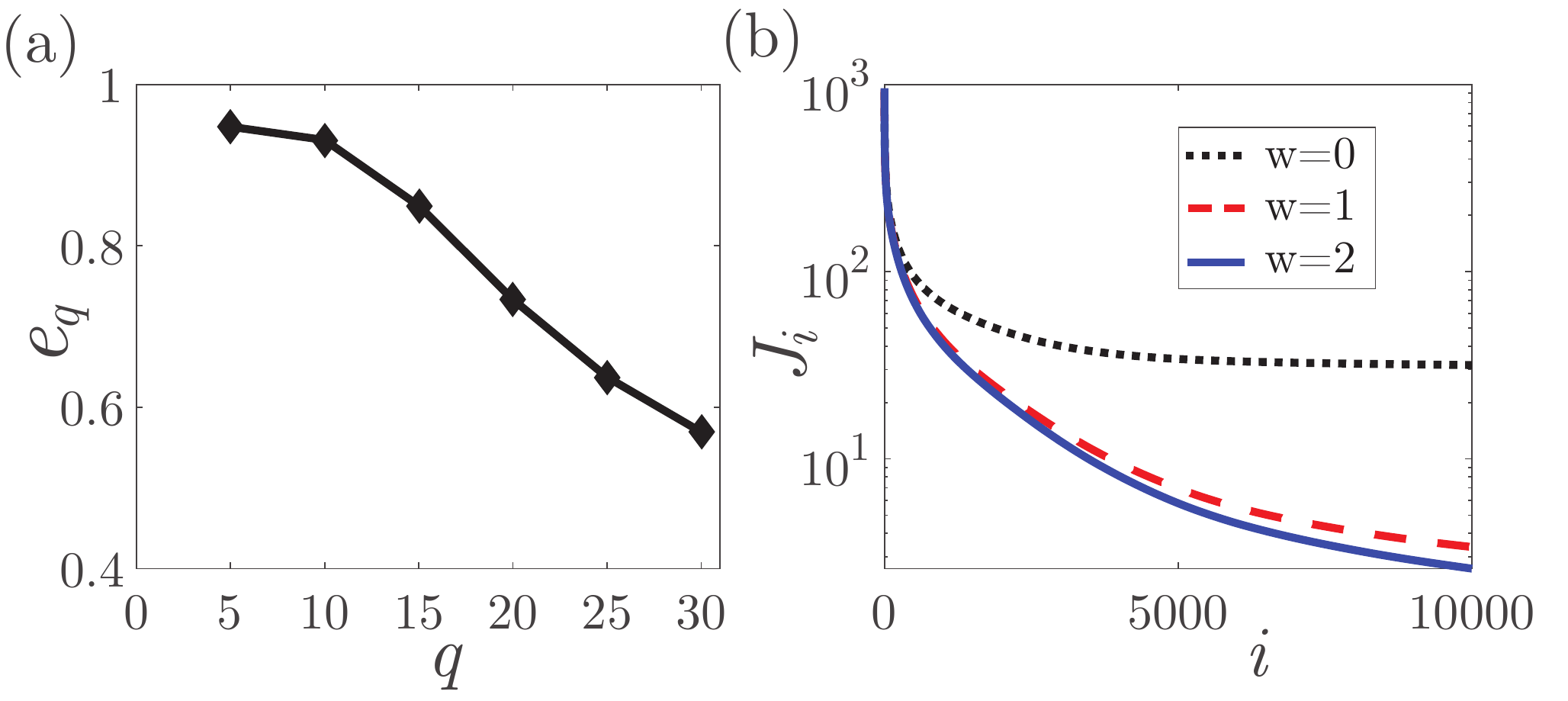}
\caption{(a) The accuracy of approximating the solution of the Lyapunov equation using \eqref{finalApproximationFaber}, as a function of $q$. (b) The convergence of the gradient projection method for several \textit{a priori} patterns. The residual $J$ is defined in \eqref{stepSize1}.}
\label{fig:AccuracyConvergence}
\end{figure}

Next, we test the performance of the second method. We investigate the convergence of  the initial guess in Method~\ref{secondMethod} as a function of the parameter $q$ in \eqref{finalApproximationFaber}. We solve the optimization problem \eqref{optimizationM} to approximate the matrix $E^{-1}$. The \textit{a priori} pattern of $E^{-1}$ is chosen as $I+E+E^2+E^3$. Such an \textit{a priori} pattern has $3.1\%$ of non-zero entries. The approximation error is $\left\|I-E\mathcal{E} \right\|_{F}=0.73$. This is a relatively crude approximation, and consequently, we are not expecting a high overall accuracy of the initial guess. Next, we compute the matrix $\mathcal{A}_{1}$.  For the sparsified Chebyshev polynomials we use an \textit{a priori} pattern $Q_{1}=I+\mathcal{A}_{1}$. The percentage of nonzero entries of $Q_{1}$ is $4.7\%$. The order of the Faber expansion is $p=30$. The fill-in of the final initial guess $\mathcal{X}_{3}$, defined by \eqref{finalApproximationFaber}, is $29\%$.  We measure the accuracy of the initial guess by \eqref{eLyapunov}, where $\hat{Z}_{k}$ is substituted by $\mathcal{X}_{3}$. The accuracy of the initial guess is shown in Fig.~\ref{fig:AccuracyConvergence}(a) as a function of the parameter $q$. We initialize the gradient projection method with such an initial guess. The convergence of the gradient projection method for several  \textit{a priori} patterns is shown in Fig.~\ref{fig:AccuracyConvergence}(b). The final approximation errors for $w=0,1,2$ are, respectively, $e=0.32, 0.26,0.25$. 

\begin{figure}[H]
\centering
 \includegraphics[scale=0.4,trim=0mm 0mm 0mm 0mm ,clip=true]{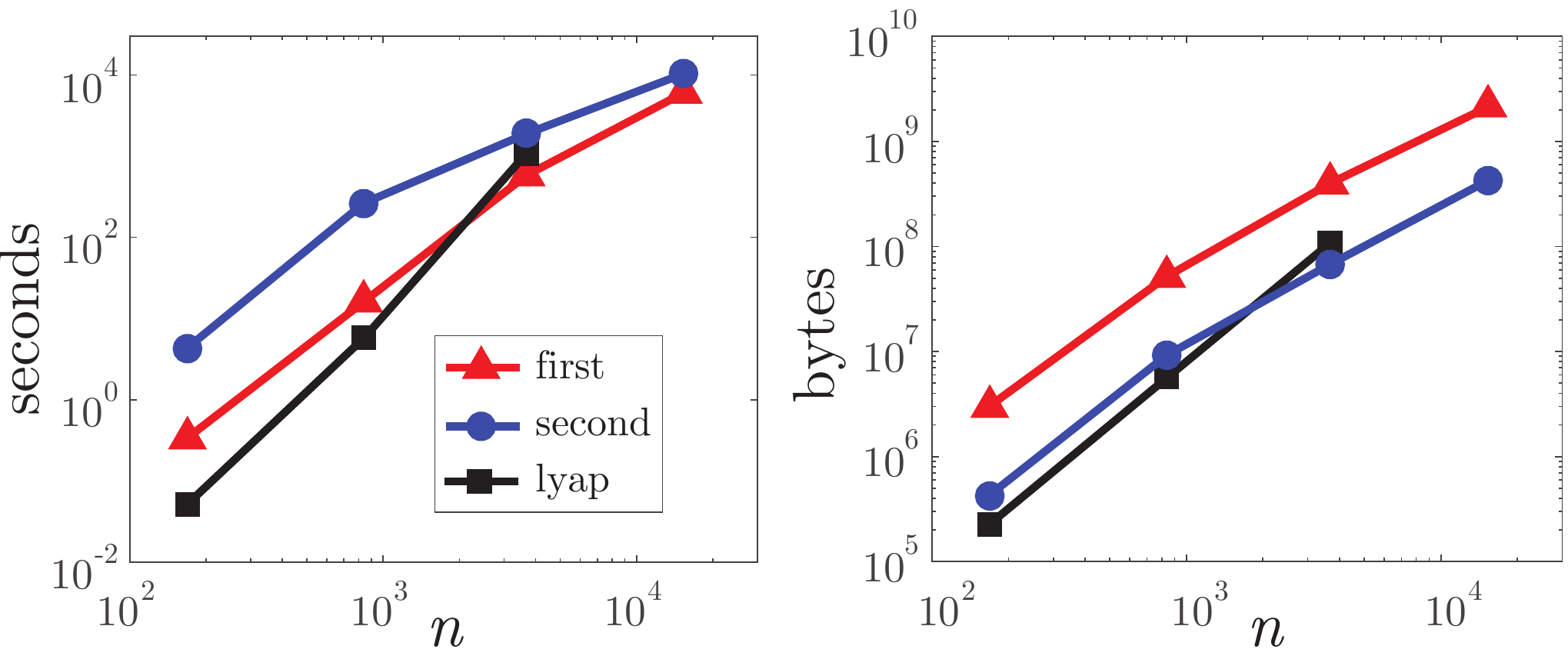}
\caption{(a) The computational and (b) memory complexities of the developed methods and the comparison with the built-in MATLAB solver $\text{lyap}\left(\cdot\right)$.}
\label{fig:Complexity}
\end{figure}

Finally, we provide some insights into the computational and memory complexities of the proposed methods. We compare the complexity of the developed methods to the complexity of the built-in MATLAB solver $\text{lyap}\left(\cdot\right)$.  The computational and memory complexities of the proposed solvers for the GL equation are shown in Fig.~\ref{fig:Complexity} for $n=168,841,3687,15379$. The exact solution obtained using the $\text{lyap}\left(\cdot\right)$ function can only be computed for $n\le 3687$. For $n=15379$, the $\text{lyap}\left(\cdot\right)$ function cannot be used because the computer runs out of RAM, and even if the computer had more RAM, it would take a significant amount of time to compute such a solution. For the first method, the residual convergence tolerance of the CGLS method is $10^{-5}$. We choose $w=1$. For these parameters, the accuracies of the converged solutions for $n=168,841,3687$ are, respectively, $4.4\cdot 10^{-4}, 0.02, 0.28$, whereas the percentages of the non-zero entries for $n=168,841,3687,15379$ are, respectively $24.5, 11.3, 3.9, 1.2$. The accuracy for $n=3687$ degrades due to the fact that the \textit{a priori} pattern generated for $w=1$ cannot capture the dominant entries as accurately as in the cases of $n=168$ and $n= 841$. This is most likely due to the increased condition number of the matrix $M$ (which can only be computed for $n=168$). To improve the accuracy, \textit{a priori} patterns $w\ge 1$ need to be used. This will increase the complexity shown in Fig.~\ref{fig:Complexity}. However, this increase will not be significant for $w=2$ and $w=3$, and the general growth trend of the complexity shown in the figure will be preserved (for larger values of $n$).
To generate the initial guess for the second method, we use the following parameters $p=20$ and $q=40$. The gradient projection method is computed for $w=1$ and computations are stopped when the iteration number exceeds $4000$. For such parameters the accuracies of the converged solutions for $n=168,841,3687$ are, respectively, $e=0.16, 0.15, 0.8$. The accuracy of the final solution can be improved by increasing both $p$ and $q$ together with the number of iterations of the gradient projection method. In our future research, we will investigate the optimal selection of these parameters that should produce high accuracy, while at the same time preserving the computational efficiency. From Fig.~\ref{fig:Complexity} we can observe that both developed approaches have lower computational and memory complexities than the MATLAB $\text{lyap}(\cdot)$ function whose computational and memory complexities are $O(n^3)$ and $O(n^2)$, respectively. It can be observed that the memory complexity of the second approximation method is in the order of magnitude smaller than the memory complexity of the first method. This is due to the fact that in the second method all operations are performed on $n\times n$ sparse matrices. On the other hand, the first method has a lower computational complexity. Since it operates directly on sparse $n\times n$ matrices, the second method can be used for problems whose dimensions are in the order of $n=10^5$ (on a computer similar to ours that has at least $16$GB RAM).  \textit{It should be noted that in the first method we are not using the preconditioning techniques.} The preconditioning techniques can significantly reduce the condition number of $M_{1}$, and consequently, they can significantly reduce the computational complexity of the first method. Similarly, the convergence rate of the projected gradient method can be increased by using scaled gradient methods~ \cite{bertsekas1999nonlinear}. Finally, the MATLAB implementation of our methods is relatively slow, and significant computational savings can be achieved by implementing the methods in the C or C$++$ programming languages.

\section{Conclusions and Future Work}
\label{sectionConclusions}
In this paper we have presented two computationally efficient methods for computing (sparse) banded approximate solutions of the generalized Lyapunov equations with banded matrices. These methods, together with the inexact Newton method for solving the generalized Riccati equations, were used to derive a novel computationally efficient approach for optimal control of finite element or finite difference state-space models.  The proposed approach is able to compute a (sparse) banded feedback matrix of the linear-quadratic optimal controller. The future research will be oriented toward the analysis of the influence of the approximation errors on the stability and performance of the inexact Newton method for approximately solving the generalized Riccati equation. 

 The proposed control approach is based on the discretized PDE models. The finite element discretization errors  depend on the type of the basis functions, as well as on the grid density. For coarse grids, the discretization errors might degrade the performance of the optimal control algorithm. These discretization errors can be seen as model uncertainties and their effect can be ameliorated using the robust control techniques, see for example~\cite{zhou1996robust}. The future research directions are to investigate the influence of the discretization errors on the control performance and to develop a method for selecting the weighting matrices (the matrices $Q$ and $R$ in the cost function \eqref{theCostFunction}), such that the robustness of the method with respect to these errors is improved. 

While preparing the final manuscript, several computationally efficient methods for approximating the solutions of the non-generalized Lyapunov equations with sparse or low-rank matrices appeared~\cite{massei2017solving,palitta2017numerical,kressner2017low}. In our future work, we will investigate the possibility of generalizing these approaches to descriptor Lyapunov equations. Some of the ideas on which these methods have been founded can be used to additionally increase the computational efficiency and numerical stability of the methods developed in this paper. In order to develop even faster methods for generating the initial guess in Method~\ref{secondMethod}, we will further investigate the possibilities and numerical performance of the matrix function approximation algorithms~\cite{benzi2017approximation,hale2008computing,higham2008functions}. 

\section{Appendix}

For presentation clarity, let us denote the coefficient $a_{l,j}$ in \eqref{faberSeries} simply by $a_{l}$, and the coefficients $c_{2,j}$, $c_{3,j}$, and $c_{4,j}$ simply by $c_{2}$, $c_{3}$, and $c_{4}$, respectively. Then, the coefficients of the truncated Faber series \eqref{faberSeries} can be approximated by \cite{bergamaschi2003efficient}:
\begin{align}
&a_{l}=\frac{1}{W}\sum_{k=0}^{W-1} g_{k} e^{-\text{i}l2\pi \frac{k}{W}},\; g_{k}= \exp (s_{k}), \notag \\ 
& s_{k}= \notag \\ 
& \left(c_{2}+\frac{c_{3}}{4c_{2}} \right)\cos \left(2\pi\frac{k}{W}\right)+c_{4}+\text{i}\left(c_{2}-\frac{c_{3}}{4c_{2}}  \right)\sin \left( \left(2\pi\frac{k}{W}\right) \right),
\notag 
\end{align}
where $W$ is a sufficiently large positive integer and $\text{i}$ is the imaginary unit.

\section*{Acknowledgments}
This work was supported by the European Research Council (ERC) (339681) and PSC CUNY Award A (60156-00 48).

\section*{References}


\end{document}